\input amstex
\documentstyle{amsppt}

\magnification=\magstep1
\parskip=6pt

\define\ce{{\Cal E}}

\define\cg{{\Cal G}}
\define\ci{{\Cal I}}
\define\co{{\Cal O}}
\define\ct{{\Cal T}}
\define\pp{{\Bbb P}}

\define\cc{{\Bbb C}}
\define\zz{{\Bbb Z}}

\define\pe{{\Bbb P(\Cal E)}}
\define\he{{H(\Cal E)}}

\define\del{\delta}
\define\vs{\varSigma}

\define\cok{\operatorname{Coker}}
\define\ext{\operatorname{Ext}}
\define\len{\operatorname{length}}
\define\pic{\operatorname{Pic}}

\define\supp{\operatorname{Supp}}

\topmatter

\title 
  Rank-2 ample vector bundles on some smooth rational surfaces
\endtitle

\author
  Hironobu Ishihara
\endauthor

\address
  Department of Mathematics,             
  Tokyo Institute of Technology,
  Oh-okayama, Meguro-ku, Tokyo 152, Japan
\endaddress

\email
  ishihara\@math.titech.ac.jp
\endemail

\date
  May 13, 1996
\enddate


\abstract
  Some classification results for ample vector bundles of rank 2
  on Hirzebruch surfaces, and on Del Pezzo surfaces, are obtained.
  In particular, we classify rank-2 ample vector bundles 
  with $c_2$ less than 7 on Hirzebruch surfaces,
  and with $c_2$ less than 4 on Del Pezzo surfaces. 
\endabstract

\endtopmatter

\document

\head 
 Introduction
\endhead 

In recent years ample and spanned vector bundles with small 
Chern numbers have been studied by several authors.
Among them, Lanteri-Sommese \cite{LS} proved that 
$(S, \ce)\simeq(\pp^2, \co(1))$ when $S$ is a normal surface
and $\ce$ is an ample and spanned rank-2 vector bundle
with $c_2(\ce)=1$ on $S$.
Ballico-Lanteri \cite{BL} and Noma \cite{N1} classified 
ample and spanned rank-2 vector bundles with $c_2=2$ 
on smooth surfaces.
Noma \cite{N2} extended the classification to the case
of normal Gorenstein surfaces.

Motivated by the results above, we attempt to classify
ample vector bundles with small $c_2$ on surfaces
without spannedness.
As the first step, we consider rank-2 ample vector bundles
on Hirzebruch surfaces, and on Del Pezzo surfaces,
in the present paper.
We obtain classification results for rank-2 ample vector bundles
with $c_2$ less than 7 on Hirzebruch surfaces, 
and with $c_2$ less than 4 on Del Pezzo surfaces.
Note that we do not treat all smooth rational surfaces
because of technical difficulty.

This paper is organized as follows.
In \S1 we collect some preliminary results.
In \S2 we study ample vector bundles $\ce$ of rank 2
on $e$-th Hirzebruch surfaces.
We see that $c_2(\ce)\ge e+2$, and classification results for $\ce$
with $e+2\le c_2(\ce)\le e+6$ are given.
As a corollary, we obtain a classification of $\ce$ with $c_2(\ce)\le 6$.
In \S3 we study ample vector bundles $\ce$ of rank 2
on Del Pezzo surfaces of degree $d\le 7$.
We see that $c_2(\ce)\ge d$, and classification results for $\ce$
with $c_2(\ce)=d,~d+1$ are given.
A partial classification result for $\ce$
with $c_2(\ce)=d+2$ is also given.
As a corollary, we obtain a classification of $\ce$ with $c_2(\ce)\le 3$.
In \S4 we study ample vector bundles $\ce$ of rank 2 on $\pp^2$.
We see that $c_2(\ce)\ge c_1(\ce)-1$, and classification results for $\ce$
with $c_1(\ce)-1\le c_2(\ce)\le c_1(\ce)+2$ are given.
Then, using the classification of $\ce$ with $c_1(\ce)\le 3$,
we obtain a classification of $\ce$ with $c_2(\ce)\le 6$.
 
\subhead
 Acknowledgment
\endsubhead

The author would like to express his gratitude to 
Professor Takao Fujita
for his useful suggestions and valuable comments.

\subhead 
 Notation and Terminology 
\endsubhead

Basically we follow the notation and terminology of \cite{H}.
We work over the complex number field $\cc$.
Vector bundles are identified with the locally free sheaves 
of their sections, and 
line bundles are also identified with the linear equivalence classes 
of Cartier divisors.
The tensor products of line bundles are usually denoted additively,
while we use multiplicative notation for intersection products.
The linear equivalence classes are often denoted by $[~]$.
We use $=$ (resp. $\equiv$) for linear (resp. numerical) equivalence.

A line bundle $L$ on a variety $X$ is called nef if
$LC\ge 0$ for every irreducible curve $C$ in $X$.
For a morphism $f:Y\to X$, we often denote $f^*L$ by $L_Y$, 
or sometimes by $L$, when there is no fear of confusion.
For a vector bundle $\ce$ on $X$, we denote by $\pe$
the associated projective space bundle and by $\he$
the tautological line bundle on $\pe$ in the sense of \cite{Fjt2}.
We say that $\ce$ is ample if $\he$ is ample.
The determinant $\det\ce$ of $\ce$ and the first Chern class 
$c_1(\ce)$ of $\ce$ are used interchangeably.
The canonical bundle of a smooth surface $S$ is denoted by $K_S$.
For an ample line bundle $A$ on $S$, the sectional genus $g(S, A)$
(or $g(A)$ for short) of the pair $(S, A)$ is given by the formula
$2g(S, A)-2=(K_S+A)A$.
For a closed subscheme $Z$ of $S$ with the ideal sheaf $\ci_Z$,
we set $\deg Z:=\len(\co_S/\ci_Z)$.

\head
 \S1. Preliminaries
\endhead

In this section we collect some preliminary results
that will be used frequently.

\proclaim{Theorem 1.1 {\rm (Lanteri-Palleschi \cite{LP, Remark 1.3})}}
 Let $A$ be an ample line bundle on a smooth surface $S$.  
 If $K_S+A$ is not nef, then $(S, A)$ is one of the following: 
 \roster 
 \item $(S, A)\simeq(\pp^2, \co(1))$ or $(\pp^2, \co(2))$;
 \item $S$ is a $\pp^1$-bundle over a smooth curve 
       and $A_F=\co_{\pp^1}(1)$ for every fiber $F$ of the ruling.
 \endroster
\endproclaim

For the proof of this theorem, Mori's cone theorem 
\cite{M, Theorem (1.4)} and the classification theorem
of extremal rational curves \cite{M, Theorem (2.1)} are essential.

Using these two theorems, we obtain a generalization of (1.1).

\proclaim{Proposition 1.2}
 Let $\ce$ be an ample vector bundle of rank $r$ on a smooth surface $S$.  
 If $rK_S+c_1(\ce)$ is not nef, then we have one of the following: 
 \roster 
 \item $S\simeq\pp^2$ and $c_1(\ce)=\co(a)$ ($r\le a<3r$);
 \item $S$ is a $\pp^1$-bundle over a smooth curve 
       and $r\le c_1(\ce)\cdot F<2r$.
 \endroster
\endproclaim

\demo{Proof}
Suppose that $rK_S+c_1(\ce)$ is not nef.
By the cone theorem, there exists an extremal rational curve $C$ on $S$
such that $(rK_S+c_1(\ce))\cdot C<0$.
By the classification theorem of extremal rational curves,
we have one of the following:
 \roster 
 \item "(i)" $S\simeq\pp^2$ and $C$ is a line;
 \item "(ii)" $S$ is a $\pp^1$-bundle over a smooth curve 
                   and $C$ is one of its fibers;
 \item "(iii)" $C$ is a $(-1)$-curve on $S$.
 \endroster
Then the case (iii) is excluded and the assertion follows
by the next lemma. \qed
\enddemo

\proclaim{Lemma 1.3 {\rm (see, e.g., \cite{Fjt1, (1.3)})}}
 Let $S$ and $\ce$ be as above.
 Then $c_1(\ce)\cdot C\ge r$ for every rational curve $C$ on $S$.
\endproclaim

\proclaim{Theorem 1.4 {\rm (Kleiman \cite{K, Theorem 3})}}
 Let $\ce$ be an ample vector bundle of rank $r$ on a smooth surface.  
 Then we have $0<c_2(\ce)<c_1^2(\ce)$.
\endproclaim

\proclaim{Theorem 1.5 {\rm (Ballico \cite{Ba, Theorem 0.1})}}
 Let $\ce$ be an ample vector bundle of rank 2 on a smooth surface. 
 Then we have $c_1^2(\ce)\le (c_2(\ce)+1)^2$.
\endproclaim
 
\example{Remark 1.6}
 In fact, Ballico obtained the inequality in a general setting,
 though (1.5) is enough for our use.
\endexample

The following theorem is essential for the proof of (1.5).

\proclaim{Theorem 1.7 {\rm (Bogomolov \cite{Bo}, 
          see also \cite{R, Theorem 1})}}
 Let $\ce$ be an ample vector bundle of rank 2 on a smooth surface $S$.  
 Then $c_1^2(\ce)>4c_2(\ce)$ if and only if 
 there exists an exact sequence 
 $$ 0\to L\to\ce\to\ci_Z\otimes M\to 0, $$
 with $L$ and $M$ line bundles on $S$, 
 and $Z$ a zero-dimensional subscheme of $S$
 with sheaf of ideals $\ci_Z$, such that:
 \roster 
 \item"({\rm i})" $(L-M)^2>4\deg Z$;
 \item"({\rm ii})" $(L-M)\cdot A>0$ for every ample line bundle
                   $A$ on $S$.
 \endroster
\endproclaim
 
\example{Remark 1.8}
 If $\ce$ is ample in (1.7), then we see that 
 $M$ is ample in (1.7).

 Indeed, the assertion is clear in case $Z=\emptyset$.
 In case $Z\not=\emptyset$, let $\pi: S'\to S$ be the blowing-up
 of $S$ with respect to $\ci_Z$.
 We denote by $E$ the exceptional divisor corresponding to
 the inverse image ideal sheaf $\pi^{-1}\ci_Z\cdot\co_{S'}$.
 Then we have $\deg Z=-E^2$ and an exact sequence
 $$ 0\to[\pi^*L+E]\to\pi^*\ce\to[\pi^*M-E]\to 0 $$
 that is induced by the exact sequence in (1.7).
 For each irreducible curve $C$ on $S$,
 we denote by $C'$ the strict transform of $C$ under $\pi$.
 Then we have
 $$ M\cdot C=(\pi^*M)\cdot C'\ge (\pi^*M-E)\cdot C'>0 $$
 since $\pi^*M-E$ is a quotient bundle of $\pi^*\ce$.
 We have also $M^2>\deg Z>0$ since $(\pi^*M-E)^2>0$.
 Thus we conclude that $M$ is ample.
\endexample

\example{Remark 1.9}
 Using (1.8), we consider the equality condition in (1.5).

 Suppose that $c_1^2(\ce)=(c_2(\ce)+1)^2$ in (1.5).
 If $c_1^2(\ce)\le 4c_2(\ce)$, then we get $c_2(\ce)=1$
 and $c_1^2(\ce)=4$.
 If $c_1^2(\ce)>4c_2(\ce)$, then we get an exact sequence
 $$ 0\to L\to\ce\to\ci_Z\otimes M\to 0 $$
 as in (1.7).
 Note that $c_1(\ce)=L+M$ and $c_2(\ce)=LM+\deg Z$.
 We have 
 $$ \align 
      0&=(c_2(\ce)+1)^2-c_1^2(\ce) \tag 1.9.1 \\
       &=(LM)^2-L^2M^2+(L^2-(\deg Z+1))(M^2-(\deg Z+1))
               +(L+M)^2\cdot(\deg Z).
    \endalign 
 $$
 Since $\ce$ is ample, $c_1(\ce)$ is also ample
 and then $0<(L-M)\cdot c_1(\ce)=L^2-M^2$.
 Since $M$ is ample by (1.8), we get $L^2M^2\le (LM)^2$
 from the Hodge index theorem.
 We get also $M^2>\deg Z$ from the argument in (1.8).
 Thus from (1.9.1) we infer that
 $$ (LM)^2=L^2M^2,~M^2=\deg Z+1,~\text{and}~\deg Z=0. $$
 It follows that $M^2=1$ and then $L\equiv tM$
 for some integer $t\ge 2$.
\endexample 

\head
 \S2. On Hirzebruch surfaces
\endhead

\definition{Definition 2.1}
 A smooth surface $S$ is said to be an {\it e-th Hirzebruch surface}
 if $S\simeq\vs_e:=\pp(\co_{\pp^1}\oplus\co_{\pp^1}(-e))$
 for some non-negative integer $e$.
\enddefinition

In this section we denote by $S$ an $e$-th Hirzebruch surface,
by $H$ the tautological line bundle $H(\co_{\pp^1}\oplus\co_{\pp^1}(-e))$
on $S$, and by $F$ a fiber of the ruling $\rho:S\to\pp^1$.

Let $\ce$ be an ample vector bundle of rank 2 on $S$.
Since $\pic S\simeq\zz\cdot H\oplus\zz\cdot F$,
we set $c_1(\ce)=aH+bF$ for some integers $a$ and $b$.
We have $a=c_1(\ce)\cdot F\ge 2$ and $b-ae=c_1(\ce)\cdot H\ge 2$
because of (1.3).

First we consider the relation between $c_1^2(\ce)$ and $e$.

\proclaim{Lemma 2.2}
 Let $S$ and $\ce$ be as above. Then $-K_S\cdot c_1(\ce)\ge 2e+8$.
\endproclaim

\demo{Proof}
Since $K_S=-2H-(e+2)F$, we have
$$ \align -K_S\cdot c_1(\ce) &=(2H+(e+2)F)(aH+bF) \\
                             &=2(b-ae)+(e+2)a \\
                             &\ge 2\cdot 2+(e+2)\cdot 2 \\
                             &=2e+8. \qed
   \endalign $$
\enddemo

\proclaim{Proposition 2.3}
 Let $S$ and $\ce$ be as above.
 If $2K_S+c_1(\ce)$ is nef, then $c_1^2(\ce)\ge 8e+24$.
\endproclaim

\demo{Proof}
Suppose that $2K_S+c_1(\ce)$ is nef. Since $\ce$ is ample,
$c_1(\ce)$ is also ample and hence $(2K_S+c_1(\ce))c_1(\ce)\ge 0$.
Since $2K_S+c_1(\ce)=(a-4)H+(b-2e-4)F$,
we have $a\ge 4$ and then
$$ -K_S\cdot c_1(\ce)=2(b-ae)+(e+2)a
                     \ge 2\cdot 2+(e+2)\cdot 4
                     =4e+12.
$$
Thus we get $c_1^2(\ce)\ge -2K_S\cdot c_1(\ce)\ge 8e+24$. \qed
\enddemo

\proclaim{Proposition 2.4}
 Let $S$ and $\ce$ be as above.
 If $2K_S+c_1(\ce)$ is not nef, then $c_1^2(\ce)\ge 4e+8$.
\endproclaim

\demo{Proof}
Suppose that $2K_S+c_1(\ce)$ is not nef. 
From (1.2) we obtain $c_1(\ce)\cdot F=2$ or 3.

In case $c_1(\ce)\cdot F=2$, we have $\ce|_F\simeq\co_F(1)^{\oplus2}$
since $\ce$ is ample.
Then $(\ce\otimes[-H])_F\simeq\co_F^{\oplus2}$. 
Hence $\cg:=\rho_*(\ce\otimes[-H])$ is 
a locally free sheaf of rank 2 on $\pp^1$
and $\rho^*\cg\simeq\ce\otimes[-H]$.
We can set $\cg=\co(t_1)\oplus\co(t_2)$ for some integers 
$t_1$ and $t_2$ ($t_1\le t_2$).
Then $\ce\simeq[H+t_1F]\oplus[H+t_2F]$.
Note that $t_1>e$ and $t_2>e$ since $\ce$ is ample. 
We have
$$ c_1^2(\ce)=(2H+bF)^2
             =-4e+4b
             \ge -4e+4(2e+2)
             =4e+8. $$

In case $c_1(\ce)\cdot F=3$, we have 
$(\ce\otimes[-2H])_F\simeq\co_F\oplus\co_F(-1)$.
Hence $\cg:=\rho_*(\ce\otimes[-2H])$ is an invertible sheaf on $\pp^1$
and the morphism $\rho^*\cg\to\ce\otimes[-2H]$ is injective. 
We have
$$ \cok(\rho^*\cg\to\ce\otimes[-2H])=\det(\ce\otimes[-2H])-\rho^*\cg
                                    =-H+(b-t)F,
$$
where $t:=\deg\cg$. Then we get an exact sequence
$$ 0\to[2H+tF]\to\ce\to[H+(b-t)F]\to 0. $$
Note that $b-t>e$ since $\ce$ is ample. 
We have
$$ c_1^2(\ce)=(3H+bF)^2
             =-9e+6b
             \ge -9e+6(3e+2)
             =9e+12
             >4e+8. \qed
$$
\enddemo 

\proclaim{Theorem 2.5}
 Let $S$ be an $e$-th Hirzebruch surface
 and $\ce$ an ample vector bundle of rank 2 on $S$.
 Then $c_1^2(\ce)\ge 4e+8$, and equality holds if and only if
 $\ce\simeq[H+(e+1)F]^{\oplus2}$,
 where $H$ is the tautological line bundle on $S$ 
 and $F$ is a fiber of the ruling.

 Furthermore, if $4e+9\le c_1^2(\ce)\le 8e+12$,
 then $\ce\simeq[H+t_1F]\oplus[H+t_2F]$, where 
 $t_1$, $t_2\in\zz$, $e+1\le t_1\le t_2$, 
 and $t_1+t_2\le 3e+3$.
\endproclaim

\demo{Proof}
From (2.3) and (2.4) we obtain $c_1^2(\ce)\ge 4e+8$ immediately.
Suppose that $c_1^2(\ce)\le 8e+12$. 
Then $2K_S+c_1(\ce)$ is not nef by (2.3).
In view of the argument in (2.4), 
there are the following two possibilities:
\roster
\item"({i})" $c_1(\ce)\cdot F=2$, $\ce\simeq[H+t_1F]\oplus[H+t_2F]$
             ($t_1\le t_2$), $t_1>e$, $t_2>e$, 
              and $c_1^2(\ce)=-4e+4(t_1+t_2)$;
\item"({ii})" $c_1(\ce)\cdot F=3$, $0\to[2H+tF]\to\ce\to[H+(b-t)F]\to 0$
              is exact, $b-3e\ge 2$, $b-t>e$, and $c_1^2(\ce)=-9e+6b$.
\endroster

In the case (i) we see that $t_1+t_2\le 3e+3$ 
since $c_1^2(\ce)\le 8e+12$.

In the case (ii) we see that $e=0$ 
since $9e+12\le c_1^2(\ce)\le 8e+12$.
Then we have $b=2$ and $c_1(\ce)=3H+2F$.
Note that the condition $c_1(\ce)=3H+2F$ is equivalent to
the condition $c_1(\ce)=2H+3F$.
Hence we infer that $\ce\simeq[H+F]\oplus[H+2F]$
from the argument above. \qed
\enddemo

Using this theorem, we can classify rank-2 ample vector bundles
with small $c_1^2$ on Hirzebruch surfaces.

\proclaim{Corollary 2.6}
 Let $S$ be an $e$-th Hirzebruch surface.
 Then rank-2 ample vector bundles $\ce$ with $c_1^2(\ce)\le 16$
 on $S$ are the following:
 \roster
 \item"({i})" $c_1^2(\ce)=8$, $e=0$, and $\ce\simeq[H+F]^{\oplus2}$;
 \item"({ii})" $c_1^2(\ce)=12$, $e=0$, and $\ce\simeq[H+F]\oplus[H+2F]$;
 \item"({iii})" $c_1^2(\ce)=12$, $e=1$, and $\ce\simeq[H+2F]^{\oplus2}$;
 \item"({iv})" $c_1^2(\ce)=16$, $e=0$, and $\ce\simeq[H+F]\oplus[H+3F]$;
 \item"({v})" $c_1^2(\ce)=16$, $e=0$, and $\ce\simeq[H+2F]^{\oplus2}$;
 \item"({vi})" $c_1^2(\ce)=16$, $e=1$, and $\ce\simeq[H+2F]\oplus[H+3F]$;
 \item"({vii})" $c_1^2(\ce)=16$, $e=2$, and $\ce\simeq[H+3F]^{\oplus2}$.
 \endroster
\endproclaim

\demo{Proof}
Suppose that $c_1^2(\ce)\le 16$.
From (2.5) we get $e\le 2$;
moreover, $\ce$ is a vector bundle of the type (vii) in case $e=2$,
and $\ce$ is of the type (iii) or (vi) in case $e=1$.

In case $e=0$, $\ce$ is of the type (i) or (ii) if $c_1^2(\ce)\le 12$.
In case $e=0$ and $13\le c_1^2(\ce)\le 16$,  
$2K_S+c_1(\ce)$ is not nef by (2.3). 
From the proof of (2.4), we infer that $\ce$ is of the type (iv) or (v).
\qed
\enddemo

Next we consider the relation between $c_2(\ce)$ and $e$.
 
\proclaim{Proposition 2.7}
 Let $S$ be an $e$-th Hirzebruch surface 
 and $\ce$ an ample vector bundle of rank 2 on $S$.  
 If $c_1^2(\ce)>4c_2(\ce)$, then $c_2(\ce)\ge e+4$.
\endproclaim

\demo{Proof}
Suppose that $c_1^2(\ce)>4c_2(\ce)$.
From (1.7) we obtain an exact sequence
$$ 0\to L\to\ce\to\ci_Z\otimes M\to 0, $$
where $L, M\in\pic S$ and $\ci_Z$ is the ideal sheaf of a 
zero-dimensional subscheme $Z$ of $S$.
Note that $M$ is ample (cf. (1.8)).

(2.7.1) If $K_S+M$ is nef, we have $(-K_S)\cdot M\le M^2<LM$
since $(K_S+M)M\ge 0$ and $(L-M)M>0$.
We have also $(-K_S)\cdot L\le LM$
since $(K_S+M)L=(K_S+M)(L-M)+(K_S+M)M\ge 0$.
Hence we obtain
$$ 2e+8\le (-K_S)\cdot c_1(\ce)=(-K_S)\cdot L+(-K_S)\cdot M
                               <2LM\le 2c_2(\ce) 
$$
by (2.2). It follows that $c_2(\ce)>e+4$.
 
(2.7.2) If $K_S+M$ is not nef, we infer that $MF=1$ from (1.1)
since $M$ is ample.
Then we can set $M=H+tF$ and $L=(a-1)H+(b-t)F$
for integers $a$, $b$, and $t$. 
Note that $t>e$ since $M$ is ample.
We have
$$ LM-(e+4)=((a-1)H+(b-t)F)(H+tF)-(e+4)
           =(b-ae-2)+(a-2)t-2.
$$
Hence we see that $LM\ge e+4$ if $a\ge 3$ and $t\ge 2$.

In case $a=2$, we have $L-M=(b-2t)F$ and hence $b\ge 2t+1\ge 2e+3$.
Thus we see that $LM\ge e+4$ unless $b=2e+3$.
If $b=2e+3$, then $t=e+1$ and
we find that $c_2(\ce)=e+3+\deg Z$ and $c_1^2(\ce)=4e+12$.
This is a contradiction to the assumption $c_1^2(\ce)>4c_2(\ce)$.

In case $t=1$, we have $e=0$, and then $LM-(e+4)=a+b-6$.
This is non-negative because 
$$ 0<c_1^2(\ce)-4c_2(\ce)=2(a-2)(b-2)-4\deg Z. $$ 
 
As a result, we have $LM\ge e+4$ if $K_S+M$ is not nef.
It follows that $c_2(\ce)\ge e+4$. \qed
\enddemo

\proclaim{Theorem 2.8}
 Let $S$ be an $e$-th Hirzebruch surface 
 and $\ce$ an ample vector bundle of rank 2 on $S$.
 Then $c_2(\ce)\ge e+2$, and equality holds if and only if
 $\ce\simeq[H+(e+1)F]^{\oplus2}$,
 where $H$ is the tautological line bundle on $S$,
 and $F$ is a fiber of the ruling.
 
 Furthermore, $c_2(\ce)=e+3$ if and only if 
 $\ce\simeq[H+(e+1)F]\oplus[H+(e+2)F]$.
\endproclaim

\demo{Proof}
Assume that $c_2(\ce)\le e+3$.
From (2.7) we obtain $c_1^2(\ce)\le 4c_2(\ce)\le 4e+12$,
and hence $2K_S+c_1(\ce)$ is not nef by (2.3).
In view of the argument in (2.4),
there are the following two possibilities: 
\roster
\item"({i})" $c_1(\ce)\cdot F=2$, $\ce\simeq[H+t_1F]\oplus[H+t_2F]$
             ($t_1\le t_2$), $t_1>e$, $t_2>e$, 
              and $c_2(\ce)=t_1+t_2-e$;
\item"({ii})" $c_1(\ce)\cdot F=3$, $0\to[2H+tF]\to\ce\to[H+(b-t)F]\to 0$
              is exact, $b-3e\ge 2$, $b-t>e$, and $c_2(\ce)=2b-t-2e$.
\endroster

In the case (i) we see that $(t_1, t_2)=(e+1, e+1)$ or $(e+1, e+2)$.
Hence we have either 
$\ce\simeq[H+(e+1)F]^{\oplus2}$ and $c_2(\ce)=e+2$, or
$\ce\simeq[H+(e+1)F]\oplus[H+(e+2)F]$ and $c_2(\ce)=e+3$.
In the case (ii) we see that 
$$ e+3\ge c_2(\ce)=2b-t-2e>b-e\ge 2e+2, $$
and hence $e=0$, $c_2(\ce)=3$, $b=2$, and $t=1$.
Then we get an exact sequence
$$ 0\to [2H+F]\to\ce\to [H+F]\to 0. $$
Since $\ext^1([H+F], [2H+F])\simeq H^1(S, H)=0$,
we have $\ce\simeq[2H+F]\oplus[H+F]$.
Hence we obtain that $\ce\simeq[H+F]\oplus[H+2F]$.
\qed
\enddemo

\proclaim{Theorem 2.9}
 Let $S$, $\ce$, $H$, and $F$ be as in (2.8).

 {\rm (I)} $c_2(\ce)=e+4$ if and only if $\ce$ is one of the following:
 
 \widestnumber\item{(III-viii)}
 \roster
 \item"({I-i})" $\ce\simeq[H+(e+1)F]\oplus[H+(e+3)F]$ or  
                $[H+(e+2)F]^{\oplus2}$;
 \item"({I-ii})" $e=0$ and $\ce\simeq[H+F]\oplus[2H+2F]$;
 \item"({I-iii})" $e=1$ and $\ce\simeq[H+2F]\oplus[2H+3F]$.
 \endroster
 
 {\rm (II)} $c_2(\ce)=e+5$ if and only if $\ce$ is one of the following:
 
 \widestnumber\item{(III-viii)}
 \roster
 \item"({II-i})" $\ce\simeq[H+(e+1)F]\oplus[H+(e+4)F]$ or  
                 $[H+(e+2)F]\oplus[H+(e+3)F]$;
 \item"({II-ii})" $e=0$ and $\ce\simeq[H+F]\oplus[2H+3F]$;
 \item"({II-iii})" $e=0$ and $\ce\simeq[H+2F]\oplus[2H+F]$;
 \item"({II-iv})" $e=1$ and $\ce\simeq[H+2F]\oplus[2H+4F]$;
 \item"({II-v})" $e=1$ and there is a non-split exact sequence
                 \newline
                 $ 0\to [2H+2F]\to\ce\to [H+3F]\to 0 $;
 \item"({II-vi})" $e=2$ and $\ce\simeq[H+3F]\oplus[2H+5F]$.
 \endroster
 
 {\rm (III)} $c_2(\ce)=e+6$ if and only if $\ce$ is one of the following:
 
 \widestnumber\item{(III-viii)}
 \roster
 \item"({III-i})" $\ce\simeq[H+(e+1)F]\oplus[H+(e+5)F]$,  
                  $[H+(e+2)F]\oplus[H+(e+4)F]$, 
                  \newline 
                  or $[H+(e+3)F]^{\oplus2}$;
 \item"({III-ii})" $e=0$ and $\ce\simeq[H+F]\oplus[2H+4F]$;
 \item"({III-iii})" $e=0$ and $\ce\simeq[H+2F]\oplus[2H+2F]$;
 \item"({III-iv})" $e=0$ and $\ce\simeq[H+F]\oplus[3H+3F]$;
 \item"({III-v})" $e=0$ and there is a non-split exact sequence
                  \newline
                  $ 0\to [2H]\to\ce\to [H+3F]\to 0 $;
 \item"({III-vi})" $e=1$ and $\ce\simeq[H+2F]\oplus[2H+5F]$;
 \item"({III-vii})" $e=1$ and $\ce\simeq[H+3F]\oplus[2H+3F]$;
 \item"({III-viii})" $e=1$ and $\ce\simeq[H+2F]\oplus[3H+4F]$;
 \item"({III-ix})" $e=1$ and there is a non-split exact sequence
                   \newline
                   $ 0\to [2H+F]\to\ce\to [H+4F]\to 0 $;
 \item"({III-x})" $e=2$ and $\ce\simeq[H+3F]\oplus[2H+6F]$;
 \item"({III-xi})" $e=2$ and there is a non-split exact sequence
                   \newline
                   $ 0\to [2H+4F]\to\ce\to [H+4F]\to 0 $;
 \item"({III-xii})" $e=3$ and $\ce\simeq[H+4F]\oplus[2H+7F]$.
 \endroster

\endproclaim

\demo{Proof}
Suppose that $e+4\le c_2(\ce)\le e+6$.
The proof is divided into two parts.

(2.9.1) If $2K_S+c_1(\ce)$ is nef, then we get 
$$ c_1^2(\ce)\ge -2K_S\cdot c_1(\ce)\ge 8e+24\ge 4c_2(\ce) $$
from the proof of (2.3).
In case $c_1^2(\ce)=4c_2(\ce)$, we have
$e=0$ and $c_1^2(\ce)=-2K_S\cdot c_1(\ce)=24$.
Then we get $(2K_S+c_1(\ce))c_1(\ce)=0$,
and hence $c_1(\ce)=-2K_S$.
It follows that $c_1^2(\ce)=4K_S^2=32$, 
which is a contradiction.
Thus we obtain $c_1^2(\ce)>4c_2(\ce)$.
We argue as in the proof of (2.7).

First we get the same exact sequence
$$ 0\to L\to\ce\to\ci_Z\otimes M\to 0 $$
as that in (2.7).
If $K_S+M$ is nef, then we have
$$ 4e+12\le -K_S\cdot c_1(\ce)<2c_2(\ce)\le 2e+12, $$
a contradiction.
Hence $K_S+M$ is not nef.
Then we can set $M=H+tF$ and $L=(a-1)H+(b-t)F$ ($a, b, t\in\zz$).
We have $t>e$ and
$$ 0\le LM-(e+4)=(b-ae-2)+(a-2)t-2. $$
Note that $a\ge 4$, otherwise we have $(2K_S+c_1(\ce))\cdot F<0$.

In case $LM=e+4$, we have $a=4$, $t=1$, $e=0$, and $b=2$.
Then $c_1^2(\ce)=(4H+2F)^2=16$, a contradiction to
$c_1^2(\ce)\ge 8e+24=24$.
In case $LM=e+5$, we have $t=1$, $e=0$, and $(a, b)=(4, 3)$ or $(5, 2)$.
Then $c_1^2(\ce)=24$ or 20, a contradiction to
$c_1^2(\ce)\ge -2K_S\cdot c_1(\ce)=28$.

Thus we obtain $LM=c_2(\ce)=e+6$ and $\deg Z=0$ since $c_2(\ce)\le e+6$.
Then we have $t\le 2$.
If $t=1$, then $e=0$ and $a=b=4$ since $2K_S+c_1(\ce)$ is nef.
Hence we get an exact sequence
$$ 0\to[3H+3F]\to\ce\to[H+F]\to 0. $$
Since $\ext^1([H+F], [3H+3F])\simeq H^1(S, 2H+2F)=0$,
we have $\ce\simeq[H+F]\oplus[3H+3F]$.
This is the case (III-iv).
If $t=2$, then $a=4$, $b=4e+2$, and $e\le 1$.
We find that $e=1$ from 
$8e+24\le c_1^2(\ce)=(4H+(4e+2)F)^2=16e+16$.
Hence we get an exact sequence
$$ 0\to[3H+4F]\to\ce\to[H+2F]\to 0. $$
Since $\ext^1([H+2F], [3H+4F])\simeq H^1(S, 2H+2F)=0$,
we have $\ce\simeq[H+2F]\oplus[3H+4F]$.
This is the case (III-viii).

(2.9.2) If $2K_S+c_1(\ce)$ is not nef, 
we argue as in the proof of (2.8). 
We have $c_1(\ce)\cdot F=2$ or 3.
If $c_1(\ce)\cdot F=2$, then we obtain 
$\ce\simeq[H+t_1F]\oplus[H+t_2F]$,
where $t_1, t_2\in\zz$, $e+1\le t_1\le t_2$, and 
$2e+4\le t_1+t_2\le 2e+6$ since $e+4\le c_2(\ce)\le e+6$.
These are the cases (I-i), (II-i), and (III-i).

If $c_1(\ce)\cdot F=3$, then there is an exact sequence
$$ 0\to[2H+tF]\to\ce\to[H+(b-t)F]\to 0 $$
with the property that
$$ e+6\ge c_2(\ce)=2b-t-2e>b-e\ge 2e+2. $$
Hence we have $e\le 3$.
In case $e=3$, we have $c_2(\ce)=9$, $b=11$, and $t=7$.
Then we get an exact sequence
$$ 0\to[2H+7F]\to\ce\to[H+4F]\to 0. $$
Since $\ext^1([H+4F], [2H+7F])\simeq H^1(S, H+3F)=0$,
we have $\ce\simeq[H+4F]\oplus[2H+7F]$.
This is the case (III-xii).

In case $e=2$, we have $b=8$ or 9.
If $b=9$, then $c_2(\ce)=8$ and $t=6$.
Hence we get an exact sequence
$$ 0\to[2H+6F]\to\ce\to[H+3F]\to 0. $$
Since $\ext^1([H+3F], [2H+6F])\simeq H^1(S, H+3F)=0$,
we have $\ce\simeq[H+3F]\oplus[2H+6F]$.
This is the case (III-x).
If $b=8$, then $(c_2(\ce), t)=(7, 5)$ or $(8, 4)$.
In the former case we obtain $\ce\simeq[H+3F]\oplus[2H+5F]$,
which is the case (II-vi).
In the latter case we obtain an exact sequence
$$ 0\to[2H+4F]\to\ce\to[H+4F]\to 0, $$
which is non-split because $2H+4F$ is not ample.
This is the case (III-xi).

In case $e=1$, we have $5\le b\le 7$.
If $b=7$, then $c_2(\ce)=7$ and $t=5$.
We obtain $\ce\simeq[H+2F]\oplus[2H+5F]$,
which is the case (III-vi).
If $b=6$, then $(c_2(\ce), t)=(6, 4)$ or $(7, 3)$.
In the former case we obtain $\ce\simeq[H+2F]\oplus[2H+4F]$,
which is the case (II-iv).
In the latter case we obtain $\ce\simeq[H+3F]\oplus[2H+3F]$,
which is the case (III-vii).
If $b=5$, then $(c_2(\ce), t)=(5, 3)$, $(6, 2)$, or $(7, 1)$.
In the first case we obtain $\ce\simeq[H+2F]\oplus[2H+3F]$,
which is the case (I-iii).
In the second case we obtain an exact sequence
$$ 0\to[2H+2F]\to\ce\to[H+3F]\to 0, $$
which is non-split because $2H+2F$ is not ample.
This is the case (II-v).
In the last case we obtain an exact sequence
$$ 0\to[2H+F]\to\ce\to[H+4F]\to 0, $$
which is non-split.
This is the case (III-ix).

In case $e=0$, we have $2\le b\le 5$.
If $b=5$, then $c_2(\ce)=6$ and $t=4$.
We obtain $\ce\simeq[H+F]\oplus[2H+4F]$,
which is the case (III-ii).
If $b=4$, then $(c_2(\ce), t)=(5, 3)$ or $(6, 2)$.
In the former case we obtain $\ce\simeq[H+F]\oplus[2H+3F]$,
which is the case (II-ii).
In the latter case we obtain $\ce\simeq[H+2F]\oplus[2H+2F]$,
which is the case (III-iii).
If $b=3$, then $(c_2(\ce), t)=(4, 2)$, $(5, 1)$, or $(6, 0)$.
In the first case we obtain $\ce\simeq[H+F]\oplus[2H+2F]$,
which is the case (I-ii).
In the second case we obtain $\ce\simeq[H+2F]\oplus[2H+F]$,
which is the case (II-iii).
In the last case we obtain a non-split exact sequence
$$ 0\to[2H]\to\ce\to[H+3F]\to 0, $$
which is the case (III-v).
If $b=2$, then $c_1(\ce)=3H+2F$.
Note that the condition $c_1(\ce)=3H+2F$ 
is equivalent to the condition $c_1(\ce)=2H+3F$.
Hence we have already treated this case.
\qed
\enddemo 

\example{Remark 2.10}
 The existence of $\ce$ in the cases (II-v) and (III-xi)
 is shown by Fujisawa \cite{Fjs, Example 3.7}.
 The existence of $\ce$ in the case (III-v) can be shown similarly. 

 The existence of $\ce$ in the case (III-ix) is shown as follows.
 Let $C_0$ be the minimal section of $\rho:\vs_1\to\pp^1$.
 We fix a non-split exact sequence 
 $$ 0\to\co_{C_0}(-1)\to\co_{C_0}(1)^{\oplus2}\to\co_{C_0}(3)\to 0 
    \tag 2.10.1
 $$
 on $C_0$.
 Since $\ext^1([C_0+4F], [2C_0+F])\simeq H^1(\vs_1, C_0-3F)$ and
 $$ h^2(\vs_1, (C_0-3F)-C_0)=h^2(\vs_1, -3F)=h^0(\vs_1, -2C_0)=0, $$
 we have a non-trivial extention
 $$ 0\to[2C_0+F]\to\ce\to[C_0+4F]\to 0 $$
 whose restriction to $C_0$ is (2.10.1).
 Then we see that $c_1^2(\ce)=21$ and $c_2(\ce)=7$. 

 We show that the tautological line bundle $\he$ on $\pe$ is ample.
 Note that $\he^3=c_1^2(\ce)-c_2(\ce)=14$.
 The surjection $\ce\to[C_0+4F]$ above determines a section $Z$
 of the projection $p:\pe\to S$.
 Then $Z\in|\he-p^*(2C_0+F)|$ and $\he|_Z\simeq C_0+4F$ is ample.
 In addition, $\he|_{p^{-1}(C_0)}$ is ample since
 $\ce|_{C_0}=\co_{C_0}(1)^{\oplus2}$.

 Let $W$ be an arbitrary irreducible surface in $\pe$
 with the property that $W\not=Z$ and $W\not=p^{-1}(C_0)$.
 Then we infer that $|\he|_W|=|[Z]_W+[p^*(2C_0+F)]_W|$
 has a non-zero member for some $F$.
 Let $C$ be an arbitrary irreducible curve in $\pe$
 with the property that $C\not\subset Z\cup p^{-1}(C_0)$.
 Then we see that $\he\cdot C=Z\cdot C+p^*(2C_0+F)\cdot C>0$.
 
 We thus conclude that $\ce$ is ample in view of the Nakai criterion.
\endexample  

Using the results above, we can classify rank-2 ample vector bundles
with small $c_2$ on Hirzebruch surfaces.

\proclaim{Corollary 2.11}
 Let $S$ be an $e$-th Hirzebruch surface.
 Then rank-2 ample vector bundles $\ce$ with $c_2(\ce)\le 6$
 on $S$ are the following:
 \roster 
 \item $c_2(\ce)=2$, $e=0$, and $\ce\simeq[H+F]^{\oplus2}$;
 \item $c_2(\ce)=3$, $e=0$, and $\ce\simeq[H+F]\oplus[H+2F]$;
 \item $c_2(\ce)=3$, $e=1$, and $\ce\simeq[H+2F]^{\oplus2}$;
 \item $c_2(\ce)=4$, $e=0$, and $\ce\simeq[H+F]\oplus[H+3F]$;
 \item $c_2(\ce)=4$, $e=0$, and $\ce\simeq[H+2F]^{\oplus2}$;
 \item $c_2(\ce)=4$, $e=0$, and $\ce\simeq[H+F]\oplus[2H+2F]$;
 \item $c_2(\ce)=4$, $e=1$, and $\ce\simeq[H+2F]\oplus[H+3F]$;
 \item $c_2(\ce)=4$, $e=2$, and $\ce\simeq[H+3F]^{\oplus2}$;
 \item $c_2(\ce)=5$, $e=0$, and $\ce\simeq[H+F]\oplus[H+4F]$;
 \item $c_2(\ce)=5$, $e=0$, and $\ce\simeq[H+2F]\oplus[H+3F]$;
 \item $c_2(\ce)=5$, $e=0$, and $\ce\simeq[H+F]\oplus[2H+3F]$;
 \item $c_2(\ce)=5$, $e=0$, and $\ce\simeq[H+2F]\oplus[2H+F]$;
 \item $c_2(\ce)=5$, $e=1$, and $\ce\simeq[H+2F]\oplus[H+4F]$;
 \item $c_2(\ce)=5$, $e=1$, and $\ce\simeq[H+3F]^{\oplus2}$;
 \item $c_2(\ce)=5$, $e=1$, and $\ce\simeq[H+2F]\oplus[2H+3F]$;
 \item $c_2(\ce)=5$, $e=2$, and $\ce\simeq[H+3F]\oplus[H+4F]$;
 \item $c_2(\ce)=5$, $e=3$, and $\ce\simeq[H+4F]^{\oplus2}$;
 \item $c_2(\ce)=6$, $e=0$, and $\ce\simeq[H+F]\oplus[H+5F]$;
 \item $c_2(\ce)=6$, $e=0$, and $\ce\simeq[H+2F]\oplus[H+4F]$;
 \item $c_2(\ce)=6$, $e=0$, and $\ce\simeq[H+3F]^{\oplus2}$;
 \item $c_2(\ce)=6$, $e=0$, and $\ce\simeq[H+F]\oplus[2H+4F]$;
 \item $c_2(\ce)=6$, $e=0$, and $\ce\simeq[H+2F]\oplus[2H+2F]$;
 \item $c_2(\ce)=6$, $e=0$, and $\ce\simeq[H+F]\oplus[3H+3F]$;
 \item $c_2(\ce)=6$, $e=0$, and there is a non-split exact sequence
                     \newline   $0\to [2H]\to\ce\to [H+3F]\to 0$;
 \item $c_2(\ce)=6$, $e=1$, and $\ce\simeq[H+2F]\oplus[H+5F]$;
 \item $c_2(\ce)=6$, $e=1$, and $\ce\simeq[H+3F]\oplus[H+4F]$;
 \item $c_2(\ce)=6$, $e=1$, and $\ce\simeq[H+2F]\oplus[2H+4F]$;
 \item $c_2(\ce)=6$, $e=1$, and there is a non-split exact sequence
                     \newline   $0\to [2H+2F]\to\ce\to [H+3F]\to 0$;
 \item $c_2(\ce)=6$, $e=2$, and $\ce\simeq[H+3F]\oplus[H+5F]$;
 \item $c_2(\ce)=6$, $e=2$, and $\ce\simeq[H+4F]^{\oplus2}$;
 \item $c_2(\ce)=6$, $e=3$, and $\ce\simeq[H+4F]\oplus[H+5F]$;
 \item $c_2(\ce)=6$, $e=4$, and $\ce\simeq[H+5F]^{\oplus2}$.
 \endroster
\endproclaim

\demo{Proof}
Suppose that $c_2(\ce)\le 6$.
Then $e+2\le c_2(\ce)\le e+6$ by (2.8).

In case $c_2(\ce)=e+2$, we have $0\le e\le 4$ and 
$\ce$ is a vector bundle  
of the type (1), (3), (8), (17), or (32) by (2.8).

In case $c_2(\ce)=e+3$, we have $0\le e\le 3$ and 
$\ce$ is of the type (2), (7), (16), or (31) by (2.8).
 
In case $c_2(\ce)=e+4$, we have $0\le e\le 2$ and 
$\ce$ is of the type 
(4), (5), (6), (13), (14), (15), (29), or (30) by (2.9). 

In case $c_2(\ce)=e+5$, we have $0\le e\le 1$ and 
$\ce$ is of the type 
(9), (10), (11), (12), (25), (26), (27), or (28) by (2.9). 

In case $c_2(\ce)=e+6$, we have $e=0$ and 
$\ce$ is of the type 
(18), (19), (20), (21), (22), (23), or (24) by (2.9). 
\qed
\enddemo

Finally we state a classification result concerning 
$c_1^2(\ce)$ and $c_2(\ce)$.

\proclaim{Corollary 2.12}
 Let $S$ be an $e$-th Hirzebruch surface and
 $\ce$ an ample vector bundle of rank 2 on $S$.
 If $\del(\ce):=(c_2(\ce)+1)^2-c_1^2(\ce)\le 16$,
 then $\ce$ is one of the following:
 \roster
 \item"({i})" $\del(\ce)=1$, $e=0$, and $\ce\simeq[H+F]^{\oplus2}$;
 \item"({ii})" $\del(\ce)=4$, $e=0$, and $\ce\simeq[H+F]\oplus[H+2F]$;
 \item"({iii})" $\del(\ce)=4$, $e=1$, and $\ce\simeq[H+2F]^{\oplus2}$;
 \item"({iv})" $\del(\ce)=7$, $e=0$, and $\ce\simeq[H+F]\oplus[2H+2F]$;
 \item"({v})" $\del(\ce)=9$, $e=0$, and $\ce\simeq[H+F]\oplus[H+3F]$;
 \item"({vi})" $\del(\ce)=9$, $e=0$, and $\ce\simeq[H+2F]^{\oplus2}$;
 \item"({vii})" $\del(\ce)=9$, $e=1$, and $\ce\simeq[H+2F]\oplus[H+3F]$;
 \item"({viii})" $\del(\ce)=9$, $e=2$, and $\ce\simeq[H+3F]^{\oplus2}$;
 \item"({ix})" $\del(\ce)=12$, $e=0$, and $\ce\simeq[H+F]\oplus[2H+3F]$;
 \item"({x})" $\del(\ce)=15$, $e=1$, and $\ce\simeq[H+2F]\oplus[2H+3F]$;
 \item"({xi})" $\del(\ce)=16$, $e=0$, and $\ce\simeq[H+F]\oplus[H+4F]$;
 \item"({xii})" $\del(\ce)=16$, $e=0$, and $\ce\simeq[H+2F]\oplus[H+3F]$;
 \item"({xiii})" $\del(\ce)=16$, $e=1$, and $\ce\simeq[H+2F]\oplus[H+4F]$;
 \item"({xiv})" $\del(\ce)=16$, $e=1$, and $\ce\simeq[H+3F]^{\oplus2}$;
 \item"({xv})" $\del(\ce)=16$, $e=2$, and $\ce\simeq[H+3F]\oplus[H+4F]$;
 \item"({xvi})" $\del(\ce)=16$, $e=3$, and $\ce\simeq[H+4F]^{\oplus2}$.
 \endroster
\endproclaim

\demo{Proof}
First we note that $\del(\ce)\ge 0$ by (1.5)
without the assumption $\del(\ce)\le 16$.

From now on, we assume that $\del(\ce)\le 16$.
If $c_1^2(\ce)\le 4c_2(\ce)$, then we have $(c_2(\ce)-1)^2\le 16$
and hence $c_2(\ce)\le 5$.
By (2.11) we obtain that $\ce$ is of the type 
(i), (ii), (iii), (v), (vi), (vii), or (viii) in case $c_2(\ce)\le 4$;
$\ce$ is of the type
(xi), (xii), (xiii), (xiv), (xv), or (xvi) in case $c_2(\ce)=5$.

If $c_1^2(\ce)>4c_2(\ce)$, then from (1.7) we get an exact sequence
$$ 0\to L\to\ce\to\ci_Z\otimes M\to 0, $$
and from (1.9) we get
$$ \align 
    16&\ge(c_2(\ce)+1)^2-c_1^2(\ce) \\
      &=(LM)^2-L^2M^2+(L^2-(\deg Z+1))(M^2-(\deg Z+1))
              +(\deg Z)\cdot c_1^2(\ce).
   \endalign
$$
If $\deg Z>0$, then we see that $c_1^2(\ce)\le 16$.
From the proof of (2.7) we get
$$ 16\ge c_1^2(\ce)\ge (c_2(\ce)+1)^2-16\ge (e+6)^2-16\ge 20, $$
a contradiction.
Hence we have $\deg Z=0$ and then
$$ 16\ge (LM)^2-L^2M^2+(L^2-1)(M^2-1). $$
Note that $L^2>M^2>0$ since $\ce$ is ample.
Then we see that $1\le M^2\le 4$.

The case $M^2=1$ leads to a contradiction.
Indeed, we set $M=xH+yF$, where $x, y\in\zz$, $x>0$, and $y>xe$.
From $1=M^2=x(2y-xe)$, we get $x=1$ and then $2y-e=1$.
It follows that $e+1=2y\ge 2(e+1)$, a contradiction.

In case $2\le M^2\le 4$, we will show that $c_2(\ce)\le 5$.
From $(L^2-1)(M^2-1)\le 16$ and $L^2>M^2$, we get
$$ 3\le L^2\le 17~\text{if}~M^2=2;~
   4\le L^2\le 9~\text{if}~M^2=3;~
   5\le L^2\le 6~\text{if}~M^2=4.
$$
Then we see that 
$$ (LM)^2\le 16+L^2M^2-(L^2-1)(M^2-1)=15+L^2+M^2\le 34, $$
and hence $c_2(\ce)=LM\le 5$.
Using (2.11), we obtain that $\ce$ is a vector bundle
of the type (iv), (ix), or (x).
\qed 
\enddemo

\head
 \S3. On Del Pezzo surfaces (of degree $\le 7$)   
\endhead

\definition{Definition 3.1}
 A smooth surface $S$ is said to be a {\it Del Pezzo surface of degree d}
 if $-K_S$ is ample and $d=(-K_S)^2$.
\enddefinition

\proclaim{Proposition 3.2 {\rm (see, e.g., \cite{D, p. 27, Th\'eor\`eme 1})}}
 A Del Pezzo surface $S$ of degree $d$ is one of the following:
 \roster
 \item"({\rm i})" $d=9$ and $S\simeq\pp^2$;
 \item"({\rm ii})" $d=8$ and $S\simeq\vs_0$ or $\vs_1$;
 \item"({\rm iii})" $1\le d\le 7$ and $S$ is isomorphic to the
                    blowing-up of $\pp^2$ at $9-d$ points,
                    no three of which lie on a line,
                    no six of which lie on a conic,
                    and for $d=1$ all eight do not lie on a conic
                    that is singular at one of them. 
 \endroster
 Conversely, every surface satisfying the condition {\rm (i)}, {\rm (ii)}, 
 or {\rm (iii)} is a Del Pezzo surface of the corresponding degree.
\endproclaim

In this section we denote by $S$ a Del Pezzo surface of degree $d\le 7$
and by $\ce$ an ample vector bundle of rank 2 on $S$.
Note that we have already studied rank-2 ample vector bundles 
on $\vs_0$ or $\vs_1$ in \S2.
Rank-2 ample vector bundles on $\pp^2$ are studied in \S4.

Sometimes we specify a blowing-up $\rho:S\to\pp^2$ at $9-d$ points  
$x_1,\dots,x_{9-d}$ and denote by $E_1,\dots,E_{9-d}$ 
the exceptional curves of $\rho$.
Then 
$\pic S\simeq\zz\cdot H\oplus\zz\cdot E_1\oplus\dots\oplus\zz\cdot E_{9-d}$.

First we consider the relation between $c_1^2(\ce)$ and $d$.

\proclaim{Proposition 3.3}
 Let $S$ be a Del Pezzo surface of degree $d\le 7$
 and $\ce$ an ample vector bundle of rank 2 on $S$.  
 Then $2K_S+c_1(\ce)$ is nef.
\endproclaim

\demo{Proof}
This assertion follows from (1.2). \qed
\enddemo

\proclaim{Corollary 3.4}
 Let $S$ and $\ce$ be as above. Then $-K_S\cdot c_1(\ce)\ge 2d$.
\endproclaim
 
\demo{Proof}
Since $-K_S$ is ample, we get $(2K_S+c_1(\ce))(-K_S)\ge 0$ by (3.3).
Then the assertion easily follows. \qed
\enddemo

\proclaim{Corollary 3.5}
 Let $S$ and $\ce$ be as above. Then $c_1^2(\ce)\ge 4d$.
\endproclaim
 
\demo{Proof}
Since $c_1(\ce)$ is ample, we get $(2K_S+c_1(\ce))c_1(\ce)\ge 0$ by (3.3).
Then we obtain $c_1^2(\ce)\ge -2K_S\cdot c_1(\ce)\ge 4d$ by (3.4). \qed
\enddemo

The following proposition will be used later.

\proclaim{Proposition 3.6}
 Let $S$ be a Del Pezzo surface of degree $d\le 7$
 and $\ce$ an ample vector bundle of rank 2 on $S$.  
 If $c_1^2(\ce)\le 4d+8$, then we have one of the following:
 \roster
 \item"({i})" $c_1(\ce)=-2K_S$;
 \item"({ii})" $d=1$ and $c_1(\ce)=-3K_S$;
 \item"({iii})" $c_1(\ce)=-2K_S+C$, where $C$ is a $0$-curve
                (i.e., $C\simeq\pp^1$ and $C^2=0$).
 \endroster
\endproclaim

\demo{Proof}
Suppose that $c_1^2(\ce)\le 4d+8$. From (3.5) we get
$$ 4d\le -2K_S\cdot c_1(\ce)\le c_1^2(\ce)\le 4d+8, $$
and hence $2d\le -K_S\cdot c_1(\ce)\le 2d+4$.
Note that $K_S\cdot c_1(\ce)+c_1^2(\ce)=2g(\det\ce)-2$ is even,
where $g(\det\ce)$ is the sectional genus of $(S, \det\ce)$.

If $-K_S\cdot c_1(\ce)=2d$, then $(2K_S+c_1(\ce))(-K_S)=0$.
Hence we get $c_1(\ce)=-2K_S$ by (3.3).
 
If $-K_S\cdot c_1(\ce)=2d+4$, then $c_1^2(\ce)=4d+8$ and 
$(2K_S+c_1(\ce))c_1(\ce)=0$.
Hence we have $c_1(\ce)=-2K_S$ and then $c_1^2(\ce)=4d$, 
which is a contradiction.

If $-K_S\cdot c_1(\ce)=2d+3$, then $c_1^2(\ce)=4d+7$.
We have $(2K_S+c_1(\ce))^2=-5$, a contradiction to (3.3).
 
If $-K_S\cdot c_1(\ce)=2d+1$, 
then $c_1^2(\ce)=4d+3$, $4d+5$, or $4d+7$.
In case $c_1^2(\ce)=4d+3$, we have $(2K_S+c_1(\ce))^2=-1$, a contradiction.
In case $c_1^2(\ce)=4d+5$ or $4d+7$, 
we note that $((2d+1)K_S+d\cdot c_1(\ce))(-K_S)=0$.
By the Hodge index theorem, we get
$$ 0\ge ((2d+1)K_S+d\cdot c_1(\ce))^2=d(d\cdot c_1^2(\ce)-(2d+1)^2). $$
It follows that $c_1^2(\ce)=4d+5$ and $d=1$.
Then we have $(3K_S+c_1(\ce))(-K_S)=0$ and $(3K_S+c_1(\ce))^2=0$.
Hence we obtain $c_1(\ce)=-3K_S$. 

If $-K_S\cdot c_1(\ce)=2d+2$, then $c_1^2(\ce)=4d+4$, $4d+6$, or $4d+8$.
In case $c_1^2(\ce)=4d+4$ or $4d+6$, we have 
$(2K_S+c_1(\ce))^2<0$, a contradiction.
Hence we get $c_1^2(\ce)=4d+8$ and $(2K_S+c_1(\ce))^2=0$.
Note that $(2K_S+c_1(\ce))-K_S$ is ample.
By the base point free theorem, there is a fibration
$\varphi:S\to W$ such that $2K_S+c_1(\ce)=\varphi^*A$
for some ample line bundle $A$ on $W$ (cf. \cite{Fjt2, (0.4.15)}).
From $(2K_S+c_1(\ce))^2=0$ and $(2K_S+c_1(\ce))K_S=-2$,
we infer that $\dim W=1$, $\deg A=1$, and
a general fiber $C$ of $\varphi$ is a 0-curve.
Hence we obtain $c_1(\ce)=-2K_S+C$. \qed
\enddemo

\example{Remark 3.7} 
We make some comments on (3.6).

In the case of (i), we can say a little more.
If $d=1$, by Fujita \cite{Fjt1, (2.8)},
$\ce\simeq [-K_S]^{\oplus2}$ or
there is a non-split exact sequence
$$ 0\to\co_{S'}(F)\to\pi^*\ce\to\pi^*(-2K_S)-\co_{S'}(F)\to 0, $$
where $\pi:S'\to S$ is the blowing-up at
(possibly infinitely near) three points $y_1$, $y_2$, $y_3$
and $F$ is the sum of exceptional curves over $\{y_i\}_{i=1}^3$.
Moreover, the existence of $\ce$ has been proved with
$\{y_i\}_{i=1}^3$ in a generic position.

The case $d\ge 2$ is yet to be studied.
For every $d\le 7$, we easily see that $\ce:=[-K_S]^{\oplus2}$
satisfies $c_1^2(\ce)=4d$.

In the case of (ii) or (iii), 
we obtain classification results if $c_2(\ce)\le d+2$
(cf. (3.9), (3.11), and (3.13)).

In the case of (iii), we find that $C\in|H-E_1|$
for some blowing-up $\rho:S\to\pp^2$.
Indeed, each singular fiber of $\varphi$ in the proof of (3.6)
is the union of two $(-1)$-curves that intersect at one point.
By contracting one $(-1)$-curve in a singular fiber,
we get a Del Pezzo surface of degree $d+1$ 
and $C$ is still a 0-curve on it.
Thus we may consider only the case $d=7$, 
and then the assertion is clear.
\endexample

Next we consider the relation between $c_2(\ce)$ and $d$.

\proclaim{Proposition 3.8}
 Let $S$ be a Del Pezzo surface of degree $d\le 7$
 and $\ce$ an ample vector bundle of rank 2 on $S$.  
 If $c_1^2(\ce)>4c_2(\ce)$, then $c_2(\ce)>d$.
\endproclaim
 
\demo{Proof}
Suppose that $c_1^2(\ce)>4c_2(\ce)$.
From (1.7) and (1.8), we obtain an exact sequence
$$ 0\to L\to\ce\to\ci_Z\otimes M\to 0, $$
where $L$, $M\in\pic S$, $M$ is ample, and $\ci_Z$ is the ideal sheaf of
a zero-dimensional subscheme $Z$ of $S$.

Then $K_S+M$ is nef by (1.1).
Hence we obtain
$$ d=(-K_S)^2\le (-K_S)\cdot M\le M^2<LM\le c_2(\ce) $$
by the fact that
$(K_S+M)(-K_S)\ge 0$, $(K_S+M)M\ge 0$, and $(L-M)M>0$. \qed
\enddemo

\proclaim{Theorem 3.9}
 Let $S$ be a Del Pezzo surface of degree $d\le 7$
 and $\ce$ an ample vector bundle of rank 2 on $S$.
 Then $c_2(\ce)\ge d$, and equality holds if and only if
 $\ce\simeq[-K_S]^{\oplus2}$.
\endproclaim

For a proof of this theorem, we need the following lemma.

\proclaim{Lemma 3.10}
 Let $S$ be as above and let $A$ be the union of all $(-1)$-curves on $S$.
 If $d\le 6$, then $A$ is connected and a member of $|-m_dK_S|$,
 where
 $$ m_1=240,~m_2=28,~m_3=9,~m_4=4,~m_5=2,~\text{and}~m_6=1. $$
 If $d=7$, then $A$ is connected and a member of $|H|$,
 where $H$ is the pullback of $\co_{\pp^2}(1)$ 
 by the blowing-up $\rho:S\to\pp^2$.
\endproclaim

\demo{Proof}
We fix a blowing-up $\rho:S\to\pp^2$ and denote by 
$\{E_i\}_{i=1}^{9-d}$ the exceptional curves of $\rho$.
The number of $(-1)$-curves on $S$ are listed in the following table
(cf., e.g., \cite{D, p. 35, Table 3}):
$$\vbox{\offinterlineskip \eightpoint
  \halign{&\vrule#&\strut\ \hfil#\hfil\ \cr
  \noalign{\hrule} 
    & type $\backslash$ $d$ && 1 && 2 && 3 && 4 && 5 && 6 && 7 & \cr
  \noalign{\vskip-1pt \hrule \vskip1pt \hrule} 
    & $(0; -1)$ && 8 && 7 && 6 && 5 && 4 && 3 && 2 & \cr
  \noalign{\hrule} 
    & $(1; 1^2)$ && 28 && 21 && 15 && 10 && 6 && 3 && 1 & \cr
  \noalign{\hrule} 
    & $(2; 1^5)$ && 56 && 21 && 6 && 1 && 0 && 0 && 0 & \cr
  \noalign{\hrule} 
    & $(3; 2, 1^6)$ && 56 && 7 && 0 && 0 && 0 && 0 && 0 & \cr
  \noalign{\hrule} 
    & $(4; 2^3, 1^5)$ && 56 && 0 && 0 && 0 && 0 && 0 && 0 & \cr
  \noalign{\hrule} 
    & $(5; 2^6, 1^2)$ && 28 && 0 && 0 && 0 && 0 && 0 && 0 & \cr
  \noalign{\hrule} 
    & $(6; 3, 2^7)$ && 8 && 0 && 0 && 0 && 0 && 0 && 0 & \cr
  \noalign{\vskip-1pt \hrule \vskip1pt \hrule} 
    & total && 240 && 56 && 27 && 16 && 10 && 6 && 3 & \cr
  \noalign{\hrule}
}}$$
where a $(-1)$-curve $C$ is said to be of the type 
$(a_0; a_1^{n_1}, a_2^{n_2},\dots)$ if
$C\in|a_0H-\sum_{k=1}^{n_1} a_1E_{i_k}-
           \sum_{l=1}^{n_2} a_2E_{j_l}-\dots|
$
($\{E_{i_k}, E_{j_l},\dots\}_{k, l,\dots}$ are all distinct).

Then the assertion can be shown by simple computation. \qed
 
\enddemo

\demo{Proof of Theorem 3.9}
We may assume $c_1^2(\ce)\le 4c_2(\ce)$ because of (3.8).
Then we obtain $c_2(\ce)\ge d$ by (3.5).
Suppose that $c_2(\ce)=d$.
Then we have $c_1^2(\ce)=4d$, and hence $c_1(\ce)=-2K_S$ by (3.6).
  
Using the Riemann-Roch theorem, we get 
$\chi(\ce\otimes K_S)=2$.
We have $H^2(\ce\otimes K_S)=0$ since $\ce$ is ample.
Thus there exists a non-zero section $s\in H^0(\ce\otimes K_S)$.

Let $(s)_0$ be the scheme of zeros of $s$.
In case $\dim (s)_0\le 0$, we have $(s)_0=\emptyset$ 
since $c_2(\ce\otimes K_S)=0$.
Then the section $s$ induces an exact sequence
$$ 0\to\co_S\overset{\cdot s}\to\longrightarrow\ce
    \otimes K_S\to\det(\ce\otimes K_S)\to 0. 
$$
Since $\ext^1(\det(\ce\otimes K_S), \co_S)\simeq H^1(S, \co_S)=0$,
we obtain $\ce\simeq [-K_S]^{\otimes2}$.
 
We will show that the case $\dim (s)_0=1$ cannot occur.

In case $\dim (s)_0=1$, we denote by $Z$ 
the one-dimensional part of $(s)_0$ as a cycle.
For every $(-1)$-curve $C$, we have 
$C\subset\supp Z$ or $C\cap\supp Z=\emptyset$
since $[\ce\otimes K_S]_C\simeq\co_C^{\oplus2}$.
If $d\le 6$, then the union $A$ of all $(-1)$-curves on $S$
is an ample connected divisor by (3.10).
Since $AZ>0$, we see that $A\subset\supp Z$.
Then the section $s$ determines a non-zero section 
$s'\in H^0(\ce\otimes K_S\otimes[-A])$.
Hence $H^0(\pe, \he+p^*[(m_d+1)K_S])\not=0$,
where $p$ is the projection $\pe\to S$.
Since $\ce$ is ample, we have 
$$ 0<\he^2(\he+p^*[(m_d+1)K_S])
    =c_1^2(\ce)-c_2(\ce)+c_1(\ce)\cdot[(m_d+1)K_S]
    =(1-2m_d)d<0,
$$
a contradiction.

In case $\dim (s)_0=1$ and $d=7$, let $\rho:S\to\pp^2$ be the 
blowing-up of $\pp^2$ at two points $x_1$ and $x_2$.
Setting $H:=\rho^*\co_{\pp^2}(1)$ and $E_i:=\rho^{-1}(x_i)$
($i=1, 2$), we have $A=E_1+E_2+C_{12}\in|H|$,
where $C_{12}$ is the strict transform of the line in $\pp^2$
passing through $x_1$ and $x_2$.
Note that $A\subset\supp Z$ or $A\cap\supp Z=\emptyset$
since $A$ is connected.
If $A\cap\supp Z=\emptyset$, then $E_1Z=0$, $E_2Z=0$, and $HZ=AZ=0$.
This is a contradiction since 
$\pic S\simeq\zz\cdot H\oplus\zz\cdot E_1\oplus\zz\cdot E_2$.
Hence $A\subset\supp Z$ and then $s$ determines a non-zero section
$s'\in H^0(\ce\otimes K_S\otimes[-A])$.
Since $H^0(\pe, \he+p^*[K_S-H])\not=0$ and
$$ \he^2(\he+p^*[K_S-H])
    =c_1^2(\ce)-c_2(\ce)+c_1(\ce)\cdot[K_S-H]
    =d-6=1,
$$
every member $D$ of $|\he+p^*[K_S-H]|$ is irreducible and reduced.
Since $p(D)=S$, we see that $\dim (s')_0\le 0$.
Furthermore, we find that $c_2(\ce\otimes[K_S-H])=1$
and hence $(s')_0$ is one reduced point $x_0$.

Let $\pi:S'\to S$ be the blowing-up of $S$ at $x_0$
and denote by $E_0$ the exceptional curve of $\pi$.
Then $\pi^*s'$ determines a non-zero section
$s''\in H^0(S', \pi^*(\ce\otimes[K_S-H])\otimes[-E_0])$
such that $(s'')_0=\emptyset$.
Hence we get an exact sequence
$$ 0\to\co_S'\overset{\cdot s''}\to\longrightarrow
    \pi^*(\ce\otimes[K_S-H])\otimes[-E_0]
    \to\det(\pi^*(\ce\otimes[K_S-H])\otimes[-E_0])\to 0,
$$
and then the exact sequence
$$ 0\to\pi^*(-K_S+H)+\co_{S'}(E_0)\to\pi^*\ce\to
       \pi^*(-K_S-H)-\co_{S'}(E_0)\to 0
$$
is induced.

Let $C_{12}'$ be the strict transform of $C_{12}$ by $\pi$.
Since $\pi|_{C_{12}'}:C_{12}'\to C_{12}$ is an isomorphism,
we see that $[\pi^*\ce]_{C_{12}'}$ is ample,
and hence $[\pi^*(-K_S-H)-\co_{S'}(E_0)]_{C_{12}'}$ is ample.
But we have
$$ (\pi^*(-K_S-H)-\co_{S'}(E_0))C_{12}'
  =(-K_S-H)C_{12}-E_0\cdot C_{12}'\le 0, 
$$
a contradiction.
Thus the case $\dim (s)_0=1$ does not occur. \qed
\enddemo 
                         
\proclaim{Theorem 3.11}
 Let $S$ be a Del Pezzo surface of degree $d\le 7$
 and $\ce$ be an ample vector bundle of rank 2 on $S$.
 If $c_2(\ce)=d+1$, then we have either
 \roster
 \item"({i})" $d=1$ and $\ce\simeq[-K_S]\oplus[-2K_S]$, or
 \item"({ii})" $2\le d\le 7$ and there is a non-split
               exact sequence
               $$ 0\to\pi^*(-K_S)+\co_{S'}(E_0)\to\pi^*\ce\to
                  \pi^*(-K_S)-\co_{S'}(E_0)\to 0,
               $$
               where $\pi:S'\to S$ is the blowing-up of $S$
               at one point $x_0$ and $E_0:=\pi^{-1}(x_0)$.
 \endroster
\endproclaim

\demo{Proof}
Suppose that $c_2(\ce)=d+1$. If $c_1^2(\ce)>4c_2(\ce)$, 
then as in the proof of (3.8), 
we obtain an exact sequence
$0\to L\to\ce\to\ci_Z\otimes M\to 0$ and
$$ d=(-K_S)^2\le (-K_S)\cdot M\le M^2<LM\le c_2(\ce)=d+1. $$
From this inequality we get $\deg Z=0$ and $M=-K_S$
since $K_S+M$ is nef.
Then $L=(-K_S)+(2K_S+c_1(\ce))$ is ample and
$\ce\simeq[-K_S]\oplus L$ by $\ext^1(-K_S, L)=0$.
We find that $c_1^2(\ce)=L^2+3d+2$.
Then we have $L^2\ge d+3$ since $c_1^2(\ce)>4c_2(\ce)$.
Thus
$$ d(d+3)\le (-K_S)^2L^2\le (-K_SL)^2=(d+1)^2, $$
and then $d=1$ and $c_2(\ce)=2$.
Hence we obtain $L=-2K_S$ and $\ce\simeq[-K_S]\oplus[-2K_S]$.

If $c_1^2(\ce)\le 4c_2(\ce)$, then we have $c_1^2(\ce)\le 4d+4$.
From (3.6) we obtain that $c_1(\ce)=-2K_S$.
Then we get $\chi(\ce\otimes K_S)=1$ from $c_2(\ce)=d+1$.
We get also $h^2(\ce\otimes K_S)=0$,
and hence there exists a non-zero section $s\in H^0(\ce\otimes K_S)$.
Then we have $H^0(\pe, \he+p^*K_S)\not=0$ and
$$ 0<\he^2(\he+p^*K_S)
    =c_1^2(\ce)-c_2(\ce)+c_1(\ce)\cdot K_S
    =d-1.
$$

We infer that $\dim (s)_0\le 0$ as in the proof of (3.9).
Since $c_2(\ce\otimes K_S)=1$, 
we see that $(s)_0$ is one reduced point $x_0$.
Let $\pi:S'\to S$ be the blowing-up of $S$ at $x_0$ 
and denote by $E_0$ the exceptional curve of $\pi$.
Then we obtain an exact sequence
$$ 0\to\pi^*(-K_S)+\co_{S'}(E_0)\to\pi^*\ce\to
    \pi^*(-K_S)-\co_{S'}(E_0)\to 0
$$
by an argument similar to that in (3.9).
This exact sequence is non-split,
otherwise we have
$$ \co_{E_0}^{\oplus2}\simeq[\pi^*\ce]_{E_0}\simeq
   [\pi^*(-K_S)+\co_{S'}(E_0)]_{E_0}\oplus
   [\pi^*(-K_S)-\co_{S'}(E_0)]_{E_0}\simeq
   \co_{E_0}(-1)\oplus\co_{E_0}(1),
$$
a contradiction.
We have thus proved the theorem. \qed
\enddemo

\example{Remark 3.12}
 The existence of $\ce$ in the case (ii) of (3.11)
 is shown by Fujisawa \cite{Fjs, Example (3.11)}.
\endexample

\proclaim{Proposition 3.13}
Let $S$ and $\ce$ be as above.
If $c_2(\ce)=d+2$, then we have one of the following:
 \roster 
 \item"({i})" $c_1(\ce)=-2K_S$;
 \item"({ii})" $\ce\simeq[-K_S]\oplus[-K_S+C]$,
               where $C$ is a $0$-curve;
 \item"({iii})" $d=1$ and $\ce\simeq[-K_S]\oplus[-2K_S+C]$,
                where $C$ is a $(-1)$-curve;
 \item"({iv})" $d=1$ and $\ce\simeq[-K_S]\oplus[-3K_S]$;
 \item"({v})" $d=2$ and $\ce\simeq[-K_S]\oplus[-2K_S]$.
 \endroster
\endproclaim
 
\demo{Proof}
Suppose that $c_2(\ce)=d+2$.
We argue as in the proof of (3.11).

(3.13.1) If $c_1^2(\ce)>4c_2(\ce)$, we get an exact sequence
$0\to L\to\ce\to\ci_Z\otimes M\to 0$ and
$$ d\le (-K_S)\cdot M\le M^2<LM\le c_2(\ce)=d+2. $$
Then we have $\deg Z\le 1$ and $(-K_S)\cdot M=M^2$
since $K_S\cdot M+M^2=2g(M)-2$ is even.
It follows that $M=-K_S$.
If $\deg Z=1$, then $c_1^2(\ce)=L^2+3d+2$ and $L^2\ge d+7$
since $c_1^2(\ce)>4c_2(\ce)$.
It follows that $d(d+7)\le dL^2\le (d+1)^2$, a contradiction.
Hence we have $\deg Z=0$. 
Then $c_1^2(\ce)=L^2+3d+4$ and $L^2\ge d+5$
since $c_1^2(\ce)>4c_2(\ce)$.
Note that $L^2-d=2g(L)$ is even.
Thus we get $d(d+6)\le dL^2\le (d+2)^2$, and then $d\le 2$.

In case $d=2$, we have $L^2=8$.
From $K_S^2L^2=(K_SL)^2=16$, we obtain $L=-2K_S$, 
and hence $\ce\simeq[-K_S]\oplus[-2K_S]$.
This is the case (v).

In case $d=1$, we have $L^2=7$ or 9.
If $L^2=9$, then $K_S^2L^2=(K_SL)^2=9$.
Hence we obtain $L=-3K_S$ and $\ce\simeq[-K_S]\oplus[-3K_S]$.
This is the case (iv). 
If $L^2=7$, we set $D:=L+2K_S$.
We find $\chi(\co_S(D))=1$, and
we get $h^2(S, D)=h^0(S, K_S-D)$ by the Serre duality.
If $h^0(S, K_S-D)>0$, then the divisor $K_S-D$ is effective
and hence $0<(-K_S)(K_S-D)=-2$, a contradiction.
Thus we have $h^2(S, D)=0$, hence $h^0(S, D)\ge 1$.
Since $(-K_S)\cdot D=1$, every member $C$ of $|D|$ is 
an irreducible reduced curve.
Furthermore, $C$ is a $(-1)$-curve because $C^2=-1$.
Thus we obtain $L=-2K_S+C$, and hence $\ce\simeq[-K_S]\oplus[-2K_S+C]$.
This is the case (iii).

(3.13.2) If $c_1^2(\ce)\le 4c_2(\ce)$, then we have $c_1^2(\ce)\le 4d+8$.
Because of (3.6), there are the following three posibilities:
$c_1(\ce)=-2K_S$; $d=1$ and $c_1(\ce)=-3K_S$; 
$c_1(\ce)=-2K_S+C$, where $C$ is a 0-curve.

The second case leads to a contradiction as below.
Assume that $d=1$ and $c_1(\ce)=-3K_S$.
We get $\chi(\ce\otimes[2K_S])=1$ by Riemann-Roch.
We get also $h^2(\ce\otimes[2K_S])=h^0(\ce^{\vee}\otimes[-K_S])
                                  =h^0(\ce\otimes[2K_S])$
by Serre duality and the fact that $\ce\simeq\ce^{\vee}\otimes\det\ce$.
(The symbol $^{\vee}$ stands for the dual.)
Thus we have $h^0(\ce\otimes[2K_S])>0$, and then 
$h^0(\pe, \he+p^*[2K_S])>0$, where $p: \pe\to S$ is the projection.
Since $\ce$ is ample, we have
$$ 0<\he^2(\he+p^*[2K_S])=c_1^2(\ce)-c_2(\ce)+c_1(\ce)\cdot(2K_S)=0,
$$
a contradiction.

In case $c_1(\ce)=-2K_S+C$ ($C$ is a 0-curve), 
we will show that $\ce\simeq[-K_S]\oplus[-K_S+C]$.
Since $\chi(\ce\otimes K_S)=3$ and $h^2(\ce\otimes K_S)=0$,
there exists a non-zero section $s\in H^0(\ce\otimes K_S)$.
If $\dim(s)_0\le 0$, since $c_2(\ce\otimes K_S)=0$,
we get an exact sequence
$$ 0\to\co_S\overset{\cdot s}\to\longrightarrow\ce
    \otimes K_S\to\det(\ce\otimes K_S)\to 0, 
$$
and then the exact sequence
$$ 0\to[-K_S]\to\ce\to[-K_S+C]\to 0 $$
is induced.
We have 
$\ext^1([-K_S+C], [-K_S])\simeq H^1(S, -C)$
and we find $\chi(\co_S(-C))=0$.
We have also $h^0(S, -C)=0$ and $h^2(S, -C)=h^0(S, K_S+C)=0$
since $-K_S$ is ample.
Hence we obtain $h^1(S, -C)=0$, and then 
$c_1(\ce)\simeq[-K_S]\oplus[-K_S+C]$.

If $\dim(s)_0=1$, then we denote by $Z$ the one-dimensional part of
$(s)_0$ as a cycle. 
We fix a blowing-up $\rho:S\to\pp^2$ for which $C\in|H-E_1|$
(see (3.7)).
\enddemo

\example{Claim}
 $Z\in|t(H-E_1)|$ for some positive integer $t$.
\endexample

\demo{Proof}
Let $j$ be an integer such that $2\le j\le 9-d$.
We denote by $C_{1j}$ the $(-1)$-curve obtained by 
the strict transform of the line in $\pp^2$
passing through $x_1$ and $x_j$.
Since $C_{1j}\in|H-E_1-E_j|$ and $c_1(\ce)\cdot C_{1j}=2$,
we have $[\ce\otimes K_S]_{C_{1j}}\simeq\co_{C_{1j}}$,
and hence $C_{1j}\cap\supp Z=\emptyset$ or $C_{1j}\subset\supp Z$.
We have also $[\ce\otimes K_S]_{E_j}\simeq\co_{E_j}$
since $c_1(\ce)\cdot E_j=2$.
Hence there are the following two possibilities:
\roster
\item"({a})" $E_j\cap\supp Z=\emptyset$ for every $j$;
\item"({b})" $E_j\subset\supp Z$ for some $j$.
\endroster 

In the case (a), each irreducible component $Z_1$ of $Z$
can be written as $[Z_1]=uH+vE_1$ for some non-negative integers
$u$ and $v$.
Since $E_j\cap C_{1j}\not=\emptyset$, we have 
$C_{1j}\not\subset\supp Z$, and hence $C_{1j}\cap\supp Z=\emptyset$.
It follows that $0=C_{1j}\cdot Z_1=u+v$ and then $Z_1\in|u(H-E_1)|$.
Thus we obtain $Z\in|t(H-E_1)|$ for some positive integer $t$.

In the case (b), we infer that $C_{1j}\subset\supp Z$ 
from the argument above.
Note that $E_j+C_{1j}\in|H-E_1|$.
Let $t_j$ be the largest integer with the property that 
the divisor $Z-t_j(E_j+C_{1j})$ is effective.
If $E_j\subset\supp(Z-t_j(E_j+C_{1j}))$, 
then we have $C_{1j}\subset\supp(Z-t_j(E_j+C_{1j})-E_j)$.
This contradicts the definition of $t_j$,
and hence we see that $E_j\not\subset\supp(Z-t_j(E_j+C_{1j}))$
for every $j$.
Thus we obtain that 
$E_j\cap\supp(Z-\sum_{k=2}^{9-d}t_k(E_k+C_{1k}))=\emptyset$
for every $j$.
Then the claim follows from an argument similar to that 
in the case (a).
\qed
\enddemo

\demo{Proof of (3.13), continued}
From the claim we infer that $s$ determines a non-zero section
$s'\in H^0(\ce\otimes[K_S-t(H-E_1)])$ satisfying $\dim(s')_0\le 0$.
Since $c_2(\ce\otimes[K_S-t(H-E_1)])=0$,
we get an exact sequence
$$ 0\to\co_S\overset{\cdot s}\to\longrightarrow\ce
    \otimes[K_S-t(H-E_1)]\to\det(\ce\otimes[K_S-t(H-E_1)])\to 0, 
$$
and then the exact sequence
$$ 0\to[-K_S+t(H-E_1)]\to\ce\to[-K_S+(1-t)(H-E_1)]\to 0 $$
is induced.
Since $\ce$ is ample, $-K_S+(1-t)(H-E_1)$ is also ample, and hence $t=1$.
Then we see that 
$$ 0\to[-K_S+H-E_1]\to\ce\to[-K_S]\to 0 $$
is exact and $\ext^1([-K_S], [-K_S+H-E_1])=0$.
Hence we obtain $\ce\simeq[-K_S]\oplus[-K_S+H-E_1]$.
\qed
\enddemo

Using the results above, we can classify rank-2 ample vector bundles
with small $c_2$ on Del Pezzo surfaces.

\proclaim{Corollary 3.14}
 Let $S$ be a Del Pezzo surface of degree $d\le 7$.
 Then rank-2 ample vector bundles $\ce$ with $c_2(\ce)\le 3$
 on $S$ are the following:
 \roster
 \item $c_2(\ce)=1$, $d=1$, and $\ce\simeq[-K_S]^{\oplus2}$;
 \item $c_2(\ce)=2$, $d=1$, and $\ce\simeq[-K_S]\oplus[-2K_S]$;
 \item $c_2(\ce)=2$, $d=2$, and $\ce\simeq[-K_S]^{\oplus2}$;
 \item $c_2(\ce)=3$, $d=1$, and there is a non-split exact sequence  
                                \newline
               $0\to\co_{S'}(F)\to\pi^*\ce\to\pi^*(-2K_S)-\co_{S'}(F)\to 0$ 
                                as in {\rm (3.7)};
 \item $c_2(\ce)=3$, $d=1$, and $\ce\simeq[-K_S]\oplus[-K_S+C]$
                                where $C$ is a $0$-curve;
 \item $c_2(\ce)=3$, $d=1$, and $\ce\simeq[-K_S]\oplus[-2K_S+C]$,
                                where $C$ is a $(-1)$-curve;
 \item $c_2(\ce)=3$, $d=1$, and $\ce\simeq[-K_S]\oplus[-3K_S]$;
 \item $c_2(\ce)=3$, $d=2$, and there is a non-split exact sequence  
                                \newline
   $0\to\pi^*(-K_S)+\co_{S'}(E_0)\to\pi^*\ce\to\pi^*(-K_S)-\co_{S'}(E_0)\to 0$ 
                                as in {\rm (3.11)};
 \item $c_2(\ce)=3$, $d=3$, and $\ce\simeq[-K_S]^{\oplus2}$.
 \endroster
\endproclaim

\demo{Proof}
Suppose that $c_2(\ce)\le 3$.
Then $d\le c_2(\ce)\le d+2$ by (3.9).

In case $c_2(\ce)=d$, we have $d\le 3$ and
$\ce$ is a vector bundle of the type (1), (3), or (9) by (3.9).

In case $c_2(\ce)=d+1$, we have $d\le 2$ and
$\ce$ is of the type (2) or (8) by (3.11).

In case $c_2(\ce)=d+2$, we have $d=1$.
By (3.13) $\ce$ is of the type (5), (6), or (7) 
unless $c_1(\ce)=-2K_S$;
if $c_1(\ce)=-2K_S$, then $\ce$ is of the type (4) in view of (3.7).
\qed
\enddemo

Finally we state a classification result concerning 
$c_1^2(\ce)$ and $c_2(\ce)$.

\proclaim{Theorem 3.15}
 Let $S$ be a Del Pezzo surface of degree $d\le 7$
 and $\ce$ be an ample vector bundle of rank 2 on $S$.
 If $\del(\ce):=(c_2(\ce)+1)^2-c_1^2(\ce)\le 6$,
 then $\ce$ is one of the following:
 \roster
 \item"({i})" $\del(\ce)=0$, $d=1$, and $\ce\simeq[-K_S]\oplus[-tK_S]$,
              where $t\ge 1$ and $t\in\zz$;
 \item"({ii})" $\del(\ce)=1$, $d=2$, and $\ce\simeq[-K_S]^{\oplus2}$;
 \item"({iii})" $\del(\ce)=2$, $d=1$, and $\ce\simeq[-K_S]\oplus[-tK_S+C]$,
                where $t\ge 2$, $t\in\zz$, and $C$ is a $(-1)$-curve;
 \item"({iv})" $\del(\ce)=4$, $d=1$, and $\ce\simeq[-K_S]\oplus[-tK_S+C]$,
               where $t\ge 1$, $t\in\zz$, and $C$ is a $0$-curve;
 \item"({v})" $\del(\ce)=4$, $d=3$, and $\ce\simeq[-K_S]^{\oplus2}$;
 \item"({vi})" $\del(\ce)=6$, $d=1$, and 
                $\ce\simeq[-K_S]\oplus[-tK_S+C+C']$,
                where $t\ge 2$, $t\in\zz$, $C$ and $C'$ are $(-1)$-curves
                such that $C\cap C'=\emptyset$.
 \endroster
\endproclaim

\demo{Proof}
First we note that $\del(\ce)\ge 0$ by (1.5)
without the assumption $\del(\ce)\le 6$.
From now on, we assume that $\del(\ce)\le 6$.

If $c_1^2(\ce)\le 4c_2(\ce)$, we have $(c_2(\ce)-1)^2\le 6$
and hence $c_2(\ce)\le 3$.

If $c_1^2(\ce)>4c_2(\ce)$, then from (1.7) we get an exact sequence
$$ 0\to L\to\ce\to\ci_Z\otimes M\to 0, $$
and from (1.9) we get
$$ \align 
     6&\ge(c_2(\ce)+1)^2-c_1^2(\ce) \\
      &=(LM)^2-L^2M^2+(L^2-(\deg Z+1))(M^2-(\deg Z+1))
              +(\deg Z)\cdot c_1^2(\ce).
   \endalign
$$
If $\deg Z>0$, we see that $c_1^2(\ce)\le 6$.
From the proof of (3.8) we get
$$ 6\ge c_1^2(\ce)\ge (c_2(\ce)+1)^2-6\ge (d+3)^2-6\ge 10, $$
a contradiction.
Hence we have $\deg Z=0$ and then
$$ 6\ge (LM)^2-L^2M^2+(L^2-1)(M^2-1). $$
Since $L^2>M^2>0$, we have $1\le M^2\le 3$, and
$$ 3\le L^2\le 7~\text{if}~M^2=2;~L^2=4~\text{if}~M^2=3. $$
Thus, unless $M^2=1$, we get
$$ (LM)^2\le 6+L^2M^2-(L^2-1)(M^2-1)=5+L^2+M^2\le 14, $$
and hence $c_2(\ce)\le 3$.

In case $M^2=1$, we have
$$ d=(-K_S)^2\le (-K_S)\cdot M\le M^2=1 $$
since $K_S+M$ is nef.
It follows that $d=1$, $M=-K_S$, $L$ is ample, 
and $\ce\simeq[-K_S]\oplus L$ as in the proof of (3.11).
We will determine the type of $L$ for each $\del(\ce)$.
Note that $\del(\ce)$ is even, because
$$ 2g(\det\ce)-2=(K_S+c_1(\ce))c_1(\ce)
                =L(L-K_S)=(LK_S)(LK_S-1)-\del(\ce).
$$
If $\del(\ce)=0$, then we get $L=-tK_S$ for some integer $t\ge 2$ by (1.9).
This is the case (i).

If $\del(\ce)=2$, 4, or 6, we set
$$ D:=\cases
      L+(c_2(\ce)-1)K_S~\text{in case}~\del(\ce)=2, \\
      L+(c_2(\ce)-2)K_S~\text{in case}~\del(\ce)=4~\text{or}~6.
      \endcases
$$
Then we get $\chi(\co_S(D))\ge 1$.
We get also $h^2(S, D)=h^0(S, K_S-D)=0$ by $(-K_S)(K_S-D)<0$.
Hence $h^0(S, D)>0$ and we can regard $D$ as an effective divisor on $S$.

In case $\del(\ce)=2$, we have $(-K_S)\cdot D=1$ and $D^2=-1$.
Hence $D$ is a $(-1)$-curve $C$ and $L=-(c_2(\ce)-1)K_S+C$.
We have $c_2(\ce)\ge 3$ since $L$ is ample.
Conversely, for every $(-1)$-curve $C$ and for every integer $t\ge 2$,
$L:=-tK_S+C$ is ample and 
$\ce:=[-K_S]\oplus L$ satisfies $\del(\ce)=2$.
This is the case (iii).

In case $\del(\ce)=4$, we have $(-K_S)\cdot D=2$ and $D^2=0$.
hence $D=C$, $2C$, or $C+C'$, 
where $C$ and $C'$ are different irreducible reduced curves.

If $D=C$, then $C$ is a $0$-curve and $L=-(c_2(\ce)-2)K_S+C$.
We have $c_2(\ce)\ge 3$ since $L$ is ample.
Conversely, for every $0$-curve $C$ and for every integer $t\ge 1$,
$L:=-tK_S+C$ is ample and 
$\ce:=[-K_S]\oplus L$ satisfies $\del(\ce)=4$.
This is the case (iv).

If $D=2C$, then $2h^1(\co_C)-2=K_S\cdot C+C^2=-1$, 
which is a contradiction.

If $D=C+C'$, then $(-K_S)\cdot C=(-K_S)\cdot C'=1$ and 
$C\cdot C'\le 1$ since $(-K_S)\cdot D=2$ and $D^2=0$.
In case $C\cdot C'=0$, we have $(C^2, C^{\prime2})=(1, -1)$ or $(-1, 1)$
since $K_S\cdot C+C^2$ is even.
Replacing $C$ with $C'$ if necessary, 
we may consider only the case $(C^2, C^{\prime2})=(1, -1)$.
Then $C'$ is a $(-1)$-curve and 
$(-K_S)^2C^2=(-K_S\cdot C)^2=1$.
Hence we get $C\in|-K_S|$ and then
$0=C\cdot C'=-K_S\cdot C'=1$, a contradiction.
Thus we have $C\cdot C'=1$.
It follows that $C^2=C^{\prime2}=-1$,
and hence $C$ and $C'$ are $(-1)$-curves.
Then $D=C+C'$ is nef and $D-K_S$ is ample.
Thus, as in the proof of (3.6),
we get a $\pp^1$-fibration and we see that $D$ is a singular fiber of it.
Hence $D$ is linearly equivalent to a 0-curve,
and this case leads to the case (iv). 

In case $\del(\ce)=6$, we have $(-K_S)\cdot D=2$ and $D^2=-2$.
Similarly as above, we see that $D=C+C'$, 
where $C$ and $C'$ are different irreducible reduced curves.
Then $C^2=C^{\prime2}=-1$ and $C\cdot C'=0$. 
Hence $C$ and $C'$ are $(-1)$-curves and $L=-(c_2(\ce)-2)K_S+C+C'$.
We have $c_2(\ce)\ge 4$ since $L$ is ample.
Conversely, for every $(-1)$-curves $C$ and $C'$  
and for every integer $t\ge 2$,
$L:=-tK_S+C+C'$ is ample and 
$\ce:=[-K_S]\oplus L$ satisfies $\del(\ce)=6$.
This is the case (vi).

We thus conclude that $\ce$ is a vector bundle of the type 
(i), (iii), (iv), or (vi) in case $c_1^2(\ce)>4c_2(\ce)$ and $M^2=1$.
In other cases, we have already shown that $c_2(\ce)\le 3$.
Hence, by (3.14), $\ce$ is of the type (ii) or (v)
unless $\ce$ is one of the types above.
\qed 
\enddemo

\head
 \S4. On $\pp^2$
\endhead

In this section we consider rank-2 ample vector bundles
$\ce$ on $\pp^2$.
We always denote $\co_{\pp^2}$ by $\co$ for simplicity.
Since $\pic\pp^2\simeq\zz\cdot\co(1)$,
we regard $c_1(\ce)$ as an integer.
We have $c_1(\ce)\ge 2$ because of (1.3).

The following theorem is essentially due to Van de Ven \cite{V}.

\proclaim{Theorem 4.1 {\rm (cf. \cite{V})}}
 Let $\ce$ be an ample vector bundle of rank 2 on $\pp^2$.
 If $c_1(\ce)\le 3$, then $\ce$ is one of the following:
 \roster
 \item"({\rm i})" $c_1(\ce)=2$ and $\ce\simeq\co(1)^{\oplus2}$;
 \item"({\rm ii})" $c_1(\ce)=3$ and $\ce\simeq\co(1)\oplus\co(2)$;
 \item"({\rm iii})" $c_1(\ce)=3$ and $\ce\simeq\ct_{\pp^2}$,
                    where $\ct_{\pp^2}$ is the tangent bundle of $\pp^2$.
 \endroster
\endproclaim

\demo{Proof}
We have $c_1(\ce)\ge 2$.
If $c_1(\ce)=2$ (resp. 3), then $\ce|_L\simeq\co_L(1)^{\oplus2}$
(resp. $\co_L(1)\oplus\co_L(2)$) for every line $L$ in $\pp^2$.
Hence $\ce$ is a uniform vector bundle if $c_1(\ce)\le 3$,
and then the assertion follows by \cite{V}. \qed
\enddemo

\example{Remark 4.2}
 Note that $\ce$ is not necessarily a uniform vector bundle 
 if $c_1(\ce)\ge 4$.
 In case $c_1(\ce)=4$, we have $1\le c_2(\ce)\le 15$ because of (1.4),
 and we obtain classification results if $c_2(\ce)\le 6$
 (cf. (4.7)).

 By definition, $\ce$ is said to be {\it stable} 
 (resp. {\it semistable}) if $c_1(\ce)>2t$ (resp. $c_1(\ce)\ge 2t$)
 for every invertible subsheaf $\co(t)$ of $\ce$.

 Ample and stable vector bundles on $\pp^2$ are studied by 
 Le Potier \cite{Le}.
 In particular, these bundles with $c_1=4$ and $c_2=7, 8$
 are studied in detail.
\endexample

Next we consider the relation between $c_1(\ce)$ and $c_2(\ce)$.  

\proclaim{Theorem 4.3}
 Let $\ce$ be an ample vector bundle of rank 2 on $\pp^2$.
 Then $c_2(\ce)\ge c_1(\ce)-1$, and equality holds
 if and only if $\ce\simeq\co(1)\oplus\co(t)$ 
 for some positive integer $t$.
\endproclaim

\demo{Proof} 
From (1.5) we obtain $c_2(\ce)\ge c_1(\ce)-1$ immediately.
Suppose that $c_2(\ce)=c_1(\ce)-1$.
Then we obtain $c_1^2(\ce)\ge 4c_2(\ce)$.
If $c_1^2(\ce)=4c_2(\ce)$, we get $c_1(\ce)=2$,
and hence $\ce\simeq\co(1)^{\oplus2}$ by (4.1).
If $c_1^2(\ce)>4c_2(\ce)$, then from (1.7) we obtain an exact sequence
$$ 0\to\co(l)\to\ce\to\ci_Z(m)\to 0, $$
where $l, m\in\zz$ and $\ci_Z$ is the ideal sheaf 
of a zero-dimensional subscheme $Z$ of $\pp^2$.
Note that $l>m>0$. We have
$$ 0=c_2(\ce)-(c_1(\ce)-1)=(l-1)(m-1)+\deg Z, $$ 
and hence $m=1$ and $\deg Z=0$.
Then we get an exact sequence
$$ 0\to\co(l)\to\ce\to\co(1)\to 0, $$
which splits.
Hence we obtain $\ce\simeq\co(1)\oplus\co(l)$. 
Conversely, for every positive integer $t$, we easily see that 
$\ce:=\co(1)\oplus\co(t)$ satisfies $c_2(\ce)=c_1(\ce)-1$.
\qed
\enddemo

\proclaim{Proposition 4.4}
 Let $\ce$ be as above.
 Assume that $c_2(\ce)>c_1(\ce)-1$.
 Then $c_2(\ce)\ge 2c_1(\ce)-4$, and equality holds
 if and only if $\ce\simeq\co(2)\oplus\co(t)$ 
 for some integer $t\ge 2$.
\endproclaim

\demo{Proof}
If $c_1^2(\ce)>4c_2(\ce)$, as in the proof of (4.3),
we get an exact sequence
$$ 0\to\co(l)\to\ce\to\ci_Z(m)\to 0, $$
where $l, m\in\zz$, $l>m>0$, and $m^2>\deg Z$.  
Then we have
$$ 0<c_2(\ce)-(c_1(\ce)-1)=(l-1)(m-1)+\deg Z\le (l+m)(m-1), $$ 
and hence $m\ge 2$.
Thus we see that
$$ c_2(\ce)-(2c_1(\ce)-4)=(l-2)(m-2)+\deg Z\ge 0, $$ 
and equality holds if and only if $m=2$ and $\deg Z=0$.
Under these conditions, we obtain $\ce\simeq\co(2)\oplus\co(l)$.
Conversely, for every integer $t\ge 2$, we see that 
$\ce:=\co(2)\oplus\co(t)$ satisfies $c_2(\ce)=2c_1(\ce)-4$.

If $c_1^2(\ce)\le 4c_2(\ce)$, we have
$$ 4(c_2(\ce)-(2c_1(\ce)-4))\ge (c_1(\ce)-4)^2\ge 0, $$ 
and equality holds if and only if $c_1(\ce)=4$ and $c_2(\ce)=4$.
Under these conditions, we find that
$\chi(\ce(-2))=2$, $h^2(\ce(-2))=0$, and $h^0(\ce(-3))=0$.
Hence there exists a non-zero section $s\in H^0(\ce(-2))$
such that $\dim(s)_0\le 0$.
We have $(s)_0=\emptyset$ since $c_2(\ce(-2))=0$.
Then the section $s$ induces an exact sequence
$$ 0\to\co\overset{\cdot s}\to\longrightarrow\ce(-2)
    \to\det(\ce(-2))\to 0. 
$$
Since $\det(\ce(-2))\simeq\co$,
we obtain that $\ce\simeq\co(2)^{\oplus2}$. \qed
\enddemo

\proclaim{Theorem 4.5}
 Let $\ce$ be an ample vector bundle of rank 2 on $\pp^2$.

 {\rm (I)} $c_2(\ce)=c_1(\ce)$ if and only if either

 \widestnumber\item{(III-viii)}
 \roster
 \item"({I-i})" $\ce\simeq\co(2)^{\oplus2}$, or
 \item"({I-ii})" $\ce\simeq\ct_{\pp^2}$.
 \endroster

 {\rm (II)} $c_2(\ce)=c_1(\ce)+1$ if and only if either

 \widestnumber\item{(III-viii)}
 \roster
 \item"({II-i})" $\ce\simeq\co(2)\oplus\co(3)$, or
 \item"({II-ii})" $\ce$ is semistable, but not stable,
                      and there is an exact sequence
                      \newline
                      $0\to\co(2)\to\ce\to\ci_x(2)\to 0$,
                      where $x$ is a point of $\pp^2$.
 \endroster
 
 {\rm (III)} $c_2(\ce)=c_1(\ce)+2$ if and only if 
             $\ce$ is one of the following: 

 \widestnumber\item{(III-viii)}
 \roster
 \item"({III-i})" $\ce\simeq\co(2)\oplus\co(4)$;
 \item"({III-ii})" $\ce$ is not semistable, 
                       and there is an exact sequence
                       \newline
                       $0\to\co(3)\to\ce\to\ci_x(2)\to 0$,
                       where $x$ is a point of $\pp^2$;
 \item"({III-iii})" $\ce\simeq\ct_{\pp^2}(1)$;
 \item"({III-iv})" $\ce$ is stable 
                       and there is an exact sequence
                       $0\to\co(1)\to\ce\to\ci_Z(3)\to 0$,
                       where $Z$ is a zero-dimensional
                       subscheme of $\pp^2$ with $\deg Z=3$.
 \endroster
\endproclaim

\demo{Proof}
Suppose that $c_1(\ce)\le c_2(\ce)\le c_1(\ce)+2$.
From (4.4) we get $c_2(\ce)\ge 2c_1(\ce)-4$, and hence $c_1(\ce)\le 6$.
If $c_1(\ce)=6$, then $c_2(\ce)=8$ and $c_2(\ce)=2c_1(\ce)-4$.
Hence we obtain $\ce\simeq\co(2)\oplus\co(4)$ by (4.4).
This is the case (III-i).

If $c_1(\ce)\le 3$, we obtain $\ce\simeq\ct_{\pp^2}$ by (4.1).
This is the case (I-ii).

If $c_1(\ce)=5$, then we get $6\le c_2(\ce)\le 7$.
In case $(c_1(\ce), c_2(\ce))=(5, 6)$,
we have $c_2(\ce)=2c_1(\ce)-4$,
and hence $\ce\simeq\co(2)\oplus\co(3)$ by (4.4).
This is the case (II-i).
In case $(c_1(\ce), c_2(\ce))=(5, 7)$, we find that 
$\chi(\ce(-2))=3$, $h^2(\ce(-2))=0$, and $h^0(\ce(-4))=0$.
If $h^0(\ce(-3))>0$, there exists a non-zero section 
$s\in H^0(\ce(-3))$ such that $\dim(s)_0\le 0$.
Since $c_2(\ce(-3))=1$, we see that 
$(s)_0$ is one reduced point $x$ of $\pp^2$.
Hence the section $s$ induces an exact sequence
$$ 0\to\co\overset{\cdot s}\to\longrightarrow\ce(-3)
    \to\ci_x(-1)\to 0 
$$
since $c_1(\ce(-3))=-1$.
Then tensoring with $\co(3)$ gives an exact sequence
$$ 0\to\co(3)\to\ce\to\ci_x(2)\to 0. $$
This is the case (III-ii).
If $h^0(\ce(-3))=0$, there exists a non-zero section 
$s\in H^0(\ce(-2))$ such that $\dim(s)_0\le 0$.
Since $c_2(\ce(-2))=1$, we see that 
$(s)_0$ is one reduced point $x$ of $\pp^2$.
Let $\pi:\vs_1\to\pp^2$ be the blowing-up of $\pp^2$ at $x$
and $E$ the exceptional curve of $\pi$.
Then $\pi^*s$ determines a non-zero section
$s'\in H^0(\vs_1, \pi^*(\ce(-2))\otimes[-E])$ such that 
$(s')_0=\emptyset$.
Hence we get an exact sequence
$$ 0\to\co_{\vs_1}\overset{\cdot s'}\to\longrightarrow
    \pi^*(\ce(-2))\otimes[-E]\to
    \det(\pi^*(\ce(-2))\otimes[-E])\to 0 
$$
and then the exact sequence
$$ 0\to[\pi^*\co(2)+E]\to\pi^*\ce\to[\pi^*\co(3)-E]\to 0 $$
is induced. 
This exact sequence is non-split, otherwise we have
$$ \co_E^{\oplus2}\simeq[\pi^*\ce]_E\simeq
   [\pi^*\co(2)+E]_E\oplus[\pi^*\co(3)-E]_E\simeq
   \co_E(-1)\oplus\co_E(1),
$$
a contradiction.
Since 
$\ext^1([\pi^*\co(3)-E], [\pi^*\co(2)+E])\simeq
 H^1(\vs_1, \pi^*\co(-1)+2E)\simeq\cc^1, 
$
we see that $\ce$ is unique up to isomorphism if it exists.
On the other hand, $\ct_{\pp^2}(1)$ satisfies the condition of $\ce$.
Thus we conclude that $\ce\simeq\ct_{\pp^2}(1)$.
This is the case (III-iii). 

If $c_1(\ce)=4$, then we get $4\le c_2(\ce)\le 6$.
In case $(c_1(\ce), c_2(\ce))=(4, 4)$,
we have $c_2(\ce)=2c_1(\ce)-4$,
and hence $\ce\simeq\co(2)^{\oplus2}$ by (4.4).
This is the case (I-i).
In case $(c_1(\ce), c_2(\ce))=(4, 5)$, we find that 
$\chi(\ce(-2))=1$, $h^2(\ce(-2))=0$, and $h^0(\ce(-3))=0$.
Hence there exists a non-zero section 
$s\in H^0(\ce(-2))$ such that $\dim(s)_0\le 0$.
Since $c_1(\ce(-2))=0$ and $c_2(\ce(-2))=1$, 
the section $s$ induces an exact sequence
$$ 0\to\co\overset{\cdot s}\to\longrightarrow\ce(-2)
    \to\ci_x\to 0,
$$
where $x$ is a point of $\pp^2$.
Then tensoring with $\co(2)$ gives an exact sequence
$$ 0\to\co(2)\to\ce\to\ci_x(2)\to 0. $$
Note that $\ce$ is semistable since $h^0(\ce(-3))=0$. 
This is the case (II-ii).
In case $(c_1(\ce), c_2(\ce))\mathbreak
                             =(4, 6)$, we find that 
$\chi(\ce(-1))=4$, $h^2(\ce(-1))=0$, and $h^0(\ce(-3))=0$.
If $h^0(\ce(-2))>0$, then a non-zero section 
$s\in H^0(\ce(-2))$ induces an exact sequence
$$ 0\to\co\overset{\cdot s}\to\longrightarrow\ce(-2)
    \to\ci_Z\to 0, 
$$
where $Z$ is a zero-dimensional subscheme of $\pp^2$ 
with $\deg Z=c_2(\ce(-2))=2$.
Then we have $Z=x+x'$,
where $x$ and $x'$ are two points of $\pp^2$
(not necessarily distinct).
Let $L$ be the line passing through $x$ and $x'$.
Since $0\not=s|_L\in H^0(\ce(-2)|_L)$,
we see that $\ce(-2)|_L\simeq\co_L(t)\oplus\co_L(-t)$
for some integer $t\ge 2$.
It follows that $\ce|_L\simeq\co_L(2+t)\oplus\co_L(2-t)$,
which is a contradiction since $\ce|_L$ is ample.
Thus we get $h^0(\ce(-2))=0$.
Then there exists a non-zero section 
$s\in H^0(\ce(-1))$ such that $\dim(s)_0\le 0$.
The section $s$ induces an exact sequence
$$ 0\to\co\overset{\cdot s}\to\longrightarrow\ce(-1)
    \to\ci_Z(2)\to 0, 
$$
where $Z$ is a zero-dimensional subscheme of $\pp^2$ 
with $\deg Z=c_2(\ce(-1))=3$.
Then tensoring with $\co(1)$ gives an exact sequence
$$ 0\to\co(1)\to\ce\to\ci_Z(3)\to 0. $$
Note that $\ce$ is stable since $h^0(\ce(-2))=0$. 
This is the case (III-iv). \qed
\enddemo

\example{Remark 4.6}
 We make some comments on (4.5).

 The existence of $\ce$ in the cases (II-ii) and (III-ii)
 is shown by Fujisawa \cite{Fjs, Example (3.9)};
 furthermore, the existence of $\ce$ in the case (III-iv)
 is shown by Fujisawa \cite{Fjs, Example (3.11)}.

 In the case (III-iv), by Szurek-Wi\'sniewski 
 \cite{SW, p. 298~REMARK}, 
 $\ce(-1)$ is spanned and the evaluation
 $\co^{\oplus4}\to\ce(-1)$ induces an exact sequence
 $0\to\co(-1)^{\oplus2}\to\co^{\oplus4}\to\ce(-1)\to 0$.
\endexample
 
Using the theorems above, we can classify rank-2 ample vector bundles
with small $c_2$ on $\pp^2$.

\proclaim{Corollary 4.7}
 Rank-2 ample vector bundles $\ce$ with $c_2(\ce)\le 6$ on $\pp^2$
 are the following:
 \roster
 \item $c_2(\ce)=1$ and $\ce\simeq\co(1)^{\oplus2}$;
 \item $c_2(\ce)=2$ and $\ce\simeq\co(1)\oplus\co(2)$;
 \item $c_2(\ce)=3$ and $\ce\simeq\ct_{\pp^2}$;
 \item $c_2(\ce)=3$ and $\ce\simeq\co(1)\oplus\co(3)$;
 \item $c_2(\ce)=4$ and $\ce\simeq\co(2)^{\oplus2}$;
 \item $c_2(\ce)=4$ and $\ce\simeq\co(1)\oplus\co(4)$;
 \item $c_2(\ce)=5$, $\ce$ is semistable, but not stable,
                     and there is an exact sequence
                     \newline
                     $0\to\co(2)\to\ce\to\ci_x(2)\to 0$,
                     where $\ci_x$ is the ideal sheaf of a point
                     $x\in\pp^2$;
 \item $c_2(\ce)=5$ and $\ce\simeq\co(1)\oplus\co(5)$;
 \item $c_2(\ce)=6$, $\ce$ is stable, and there is an exact sequence
                     \newline
            $0\to\co(-1)^{\oplus2}\to\co^{\oplus4}\to\ce(-1)\to 0$;
 \item $c_2(\ce)=6$ and $\ce\simeq\co(2)\oplus\co(3)$;
 \item $c_2(\ce)=6$ and $\ce\simeq\co(1)\oplus\co(6)$.
 \endroster
\endproclaim

\demo{Proof}
Suppose that $c_2(\ce)\le 6$.
In case $c_1(\ce)\le 3$, $\ce$ is a vector bundle
of the type (1), (2), or (3) by (4.1).
In case $c_1(\ce)\ge 4$, we have 
$4\le c_1(\ce)\le c_2(\ce)+1\le 7$ by (4.3).
Then we see that $c_2(\ce)\le c_1(\ce)+2$ and hence  
$\ce$ is of the type (4), (5), (6), (7), (8), (9), (10), or (11) 
by (4.3), (4.5), and (4.6).
\qed
\enddemo

We conclude with a complement to (3.15).

\proclaim{Corollary 4.8}
 Let $\ce$ be an ample vector bundle of rank 2 on $\pp^2$.
 Assume that $\del(\ce):=(c_2(\ce)+1)^2-c_1^2(\ce)\le 24$.
 Then $\ce$ is one of the following:
 \roster
 \item"({i})" $\del(\ce)=0$ and $\ce\simeq\co(1)\oplus\co(t)$,
              where $t$ is a positive integer.
 \item"({ii})" $\del(\ce)=7$ and $\ce\simeq\ct_{\pp^2}$;
 \item"({iii})" $\del(\ce)=9$ and $\ce\simeq\co(2)^{\oplus2}$;
 \item"({iv})" $\del(\ce)=20$ and there is an exact sequence
               $0\to\co(2)\to\ce\to\ci_x(2)\to 0$ as in {\rm (4.7; 7)};
 \item"({v})" $\del(\ce)=24$ and $\ce\simeq\co(2)\oplus\co(3)$.
 \endroster
\endproclaim

\demo{Proof}
First we note that $\del(\ce)\ge 0$ by (1.5).
If $\del(\ce)=0$, then we have $c_2(\ce)=c_1(\ce)-1$.
Hence $\ce$ is a vector bundle of the type (i) by (4.3).
Conversely, for every positive integer $t$, we easily see that 
$\ce:=\co(1)\oplus\co(t)$ satisfies $\del(\ce)=0$.

Assume that $0<\del(\ce)\le 24$.
Then from (4.4) we get 
$$ c_1^2(\ce)+24\ge (c_2(\ce)+1)^2\ge (2c_1(\ce)-3)^2. $$
It follows that $c_1(\ce)\le 5$ and then $c_2(\ce)\le 6$.
Hence $\ce$ is of the type (ii), (iii), (iv), or (v) by (4.7). 
\qed
\enddemo

\enddocument